




\magnification\magstep1

\def \en{\{ e_n \}_{n=1}^{\infty}}
\def \xn{\{ x_n \}_{n=1}^{\infty}}

\def \kn{\{k_n\}_{n=1}^\infty}

\def \yn{\{ y_n \}_{n=1}^{\infty}}
\def \xk{\{ x_k \}_{k=1}^{\infty}}
\def \Xn{\{ X_n \}_{n=1}^{\infty}}
\def \En{\{ E_n \}_{n=1}^{\infty}}
\def \Ek{\{ E_k \}_{k=1}^{\infty}}
\def \Fn{\{F_n\}_{n=1}^\infty}
\def \EnFn{\{E_n\oplus F_n\}_{n=1}^\infty}

\def \zn{\{z_n\}_{n=1}^{\infty}}

\def \mk{\{ m_k \}_{k=1}^{\infty}}

\def\superscript#1{^{\raise3pt\hbox{$\scriptstyle #1$}}}
\def \seq#1#2{#1_1,\dots,#1_#2}
\def \sstwo#1{_{\lower2pt\hbox{${\scriptstyle #1}$}}}
\def\ss#1{\lower3pt\hbox{${\scriptstyle #1}$}}
\def \qed {\vrule height6pt  width6pt depth0pt}

\def \e{\epsilon }
\def \a{\alpha }

\def \d{\delta}
\def \t{\tau}

\def\lam{\lambda}

\def \nm#1{\left\|#1\right\|}

\def\sm{\smallskip}
\def\ms{\medskip}
\def\bs{\bigskip}
\def\cen{\centerline}
\def\smsk{\smallskip}
\def \Proof{\noindent {\bf Proof.\ \ }}

\def \Ball {{\rm Ball\,}}
\def \Mid {{\rm Mid}}
\def \spa{{\rm span\,}}
\def \sign{{\rm sign\,}}

\def\T{{\cal T}}
\def\Tt{{\cal T}^{2}}
\def\Tp{{\cal T}^{p}}
\def\Tpk{{\cal T}^{p\ss{k}}}
\def\Tpo{{\cal T}^{p\ss{1}}}

\def \UX{X\ss{\cal U}}
\def \UXp{X'\ss{\cal U}}
\def \UY{Y\ss{\cal U}}
\def \UYp{Y'\ss{\cal U}}
\def\U{{\cal U}}
 \def\UXn{\left(X_n\right)\ss{\cal U}}

\def\UXn{\left(X_n\right)\ss{\cal U}}
\def\UEk{\left(E_k\right)\ss{\cal U}}
\def\UFk{\left(F_k\right)\ss{\cal U}}
\def \UXz {X\ss{{\cal U},0}}
\def \UYz {Y\ss{{\cal U},0}}

\def\Lpo{L_{p\ss{1}}}
\def\Lpt{L_{p\ss{2}}}
\def\Lpk{L_{p\ss{k}}}

\def\sLpk{{\cal L}_{p\ss{k}}}

\def\sLpo{{\cal L}_{p\ss{1}}}
\def\sLpt{{\cal L}_{p\ss{2}}}
\def\sLp{{\cal L}_{p}}
\def\sLt{{\cal L}_{2}}

\def\lpj{\ell_{p\ss{j}}}
\def\lpn{\ell_{p\ss{n}}}

\def\lpk{\ell_{p\ss{k}}}
\def\lpo{\ell_{p\ss{1}}}
\def\lpt{\ell_{p\ss{2}}}
\def\pk{{p\ss{k}}}
\def\po{{p\ss{1}}}
\def\pt{{p\ss{2}}}

\def\qk{{q\ss{k}}}

\def\Lqo{L_{q\ss{1}}}
\def\Lqt{L_{q\ss{2}}}
\def\lqo{\ell_{q\ss{1}}}
\def\lqt{\ell_{q\ss{2}}}

\def\dstyle{\displaystyle}

\def\sLp{{\cal L}_{p}}
\def\sLo{{\cal L}_{1}}
\def \UX{X\ss{\cal U}}
\def \UY{Y\ss{\cal U}}
\def\U{{\cal U}}
\def\sLinf{{\cal L}_{\infty}}
\def \sninf{\{s_{n,\infty }\}_{n=1}^{\infty}}
\def\Rn#1{\hbox{{\it I\kern -0.25emR}$\sp{\,{#1}}$}}

\def \eevenn{\{ e_{2n} \}_{n=1}^{\infty}}
\def \eoddn{\{ e_{2n-1} \}_{n=1}^{\infty}}
\def \span{{\rm span\,}}

\def\implies{\Rightarrow}

\pageno=0
\ifnum\pageno<1 \footline={\hfil}\else\number\pageno\fi

\bs

\centerline{\bf BANACH SPACES DETERMINED BY THEIR UNIFORM
STRUCTURES}\footnote{}{{\it Subject classification:} 46B20, 54Hxx. {\it
Keywords:} Banach spaces, Uniform homeomorphism, Lipschitz homeomorphism}

\bs
\centerline{ by William B.~Johnson\footnote{*}{
Erna and Jacob Michael Visiting Professor, The Weizmann
Institute, 1994}\footnote{${}^\dagger$}{Supported in part by NSF
DMS 93-06376}{}\footnote{${}^\ddagger$}{Supported in part by the
U.S.-Israel Binational Science Foundation}, Joram
Lindenstrauss{${}^{\ddagger}$\footnote{$^+$}{Participant, Workshop in Linear
Analysis and Probability, Texas A\&M University}}, and Gideon
Schechtman{${}^{\ddagger}$${}^{+}$}}

\bs

\bs

\centerline{\bf Dedicated to the memory of E. Gorelik}

\bs

\bs

\centerline {\bf Abstract}

\bigskip

\noindent Following results of Bourgain and Gorelik we show that the
spaces
$\ell_p$, $1<p<\infty$, as well as some related spaces have the
following uniqueness property: If $X$ is a Banach space uniformly
homeomorphic to one of these spaces then it is linearly isomorphic
to the same space. We also prove that if a $C(K)$ space is
uniformly homeomorphic to $c_0$, then it is isomorphic to $c_0$. We
show also that there are Banach spaces which are uniformly
homeomorphic to exactly $2$  isomorphically distinct spaces.

\vfill\eject

\footline={\hss\tenrm\number\pageno\hss}

\noindent{\bf 0. Introduction}

\bs

The first result in the subject we study is the
Mazur-Ulam theorem which says that an isometry
from one Banach space onto another which takes the
origin to the origin must be linear.  This result,
which is nontrivial only when the Banach spaces are not
strictly convex, means that the structure of a Banach
space as a metric space determines the linear structure
up to translation.  On the other hand, the structure of
an infinite dimensional Banach space as a topological
space gives no information about the linear structure
of the space [Kad], [Tor].  In this paper we are
concerned with equivalence relations of Banach spaces
which lie between isometry and homeomorphism, mainly
Lipschitz equivalence and uniform homeomorphism. There
exists a considerable literature on this topic (see
[Ben2] for a nice survey to about 1983). Nevertheless,
the subject is still in its infancy, as many
fundamental, basic questions remain unanswered.  What
we find fascinating is that the subject combines
topological arguments and constructions with deep
facts from the linear structure theory of Banach
spaces. The present paper is also largely concerned
with this interplay.

Early work in this subject, especially that of Ribe
[Rib1], showed that if two Banach spaces are uniformly
homeomorphic, then each is finitely crudely
representable in the other. (Recall that $Z$ is  {\sl finitely
 crudely representable\/} in $X$
 provided there is a constant $\lam$ so that each finite
dimensional subspace
$E$ of
$Z$ is
$\lam$-isomorphic to a subspace of $X$).
 In a short but imprecise way,
this means that the two spaces have the same finite
dimensional subspaces.  Since the spaces $\ell_p$ and
$L_p(0,1)$ have the same finite dimensional subspaces,
a natural question was whether they are uniformly
homeomorphic for $1\le p <\infty$, $p\not=2$. This
question was answered in the negative first for $p=1$
by Enflo [Ben2], then for $1<p<2$ by Bourgain [Bou],
and, finally,  for $2<p<\infty$ by Gorelik [Gor].  The
main new point in Gorelik's proof is a nice
topological argument using the Schauder fixed point
theorem. In section 1 we formulate what Gorelik's
approach yields as ``The Gorelik Principle".  This
principle is most conveniently applicable for getting
information about spaces which are uniformly
homeomorphic  to a  space which has an unconditional
basis with a certain  convexity or concavity
property (see e.g. Corollary 1.7).

In section 2 we  combine Bourgain's result [Bou] and
the Gorelik Principle with structural results from the
linear theory to conclude that the linear structure of
$\ell_p$,
$1<p<\infty$, is determined by its uniform structure;
that is, the only Banach spaces which are uniformly
homeomorphic to
$\ell_p$ are those which are isomorphic to $\ell_p$.
The case $p=2$ was done twenty-five years ago by Enflo
[Enf].  The main part of section 2 is devoted to
proving that some  other spaces are determined by their
uniform structure and also to investigating the
possible number of linear structures on spaces which are
``close" to
$\ell_p$ in an appropriate local sense.

In section 3 we use the Gorelik Principle to study the
uniqueness question for $c_0$, answering in the process
a question  Aharoni asked in 1974 [Aha].

Very roughly, the passage from a uniform homeomorphism
$U$ between Banach  spaces to a linear isomorphism
involves two steps:

1. Passage from a uniformly continuous $U$ to a
Lipschitz map $F$ via the ``formula"
$F(x)=\lim_{n\to\infty} n^{-1} U(nx)$.

2. Passage from a Lipschitz map $F$ to its derivative.

Of course, both steps lead to difficulties.  In step
one a major problem is that the limit does not exist in
general and thus one is forced to use ultrafilters and
ultraprooducts.  In step two the problem is again the
the existence of derivatives.  Derivatives in the sense
of Gateaux of Lipschitz functions exist under rather
mild conditions, but they usually do not suffice; on
the other hand, derivatives in the sense of Frech\'et,
which are much more useful, exist (or at least are
known to exist) only in special cases.

It is therefore natural that in section 2 ultraproducts
are used as well as a recent result [LP] on
differentiation which ensures the existence of a
derivative in  sense which is between those of Gateaux
and Frech\'et.

It is evident that for step one of the procedure above
it is important that $U$ be defined on the entire
Banach space.  In fact, it often happens that the unit
balls of spaces $X$ and $Y$ are uniformly homeomorphic
while $X$ and $Y$ are not.   The simplest example of
this phenomenon goes back to Mazur [Maz] who noted that
the  map from $L_p(\mu)$, $1< p< \infty$, to
$L_1(\mu)$ defined by $f\mapsto |f|^p\sign f$ is a
uniform homeomorphisms between the unit balls of these
spaces, while in [Lin1] and [Enf] it is  proved that the spaces
themselves are not uniformly homeomorphic (this of course also
follows from Ribe's result mentioned above).  The Mazur map was
extended recently to more general situations by Odell and
Schlumprecht [OS].  In section 4 we obtain estimates on the
modulus of continuity of these generalized Mazur maps.  The
proofs are based on complex interpolation.

In the rather technical section 5 we combine the results
of sections 2 and 4 with  known constructions to
produce examples, for each $k=1,2,3,\dots$, of
Banach spaces which admit exactly $2^k$ linear
structures.   The most easily described examples are
certain direct sums of convexifications of Tsirelson's
space (see [Tsi], [CS]).

Although our interest is mainly in the separable
setting, we present in section 6 some results for
nonseparable spaces.

As pointed out by Ribe [Rib1], Enflo's result on
$\ell_2$ [Enf] mentioned earlier follows from the fact
that $\ell_2$ is determined by its finite dimensional
subspaces.  Section 7, which is formally independent of
the rest of the paper,  is motivated by the problem of
characterizing those Banach spaces which are determined
by their finite dimensional subspaces.  We conjecture
that only $\ell_2$ has this property.  We show that any
space which is determined by its finite dimensional
subspaces must be close to $\ell_2$, and if we
replace ``determined by its finite dimensional
subspaces" by a natural somewhat stronger property then,
besides $\ell_2$, there are other spaces which enjoy this
property.

The paper ends with a section which mentions a few of
the many open problems connected to the results of
sections 1--5.

We use standard Banach space theory language and
notation, as may be found in [LT1,2] and [T-J].

The  authors thank David
Preiss for discussions which led to the
formulation of Proposition 2.8.
\ms

\cen{*\quad\quad*\quad\quad*\quad\quad*\quad\quad*}
\ms

As is clear from this introduction, our work on this
paper was motivated by Gorelik's paper [Gor].  A few
days after he submitted the final version of [Gor],
Gorelik was fatally injured by a car while he was
jogging. We all had the privilege of knowing Gorelik,
unfortunately for only a very short time.

\bs

\noindent{\bf 1. The Gorelik Principle}

\bs

A careful reading of Gorelik's paper [Gor] leads to the
formulation of the following principle:

\smsk

\proclaim The Gorelik Principle.  A uniform homeomorphism
between Banach spaces cannot take a large ball of a finite
codimensional subspace into a small neighborhood of a subspace
of infinite codimension.

\smsk

The precise formulation of the Gorelik Principle in the form we
use is:

\proclaim Theorem 1.1. Let $U$ be a  homeomorphism from a
Banach space $X$ onto a Banach space $Y$ with uniformly
continuous inverse $V$. Suppose that
$d$,
$b$ are such that there exist a finite codimensional subspace
 $X_0$ of
$X$ and an infinite codimensional subspace    $Y_0$ of    $Y$
for which
$$
U\left[ d \Ball(X_0) \right] \subset Y_0 + b \Ball(Y).
$$
Then \
$
\omega(V,2b)\ge d/4
$, \
where \  $\omega(V,t) \equiv \displaystyle
\sup\{\nm{Vy_1-Vy_2} : \nm{y_1-y_2}\le t \}$ \  is the modulus
of uniform continuity of $V$.

The proof of the Gorelik Principle is based on two lemmas
which  are minor variations of Lemmas 5 and 6 in Gorelik's paper
[Gor]. The first lemma is obvious, while the second is a simple
consequence of Brouwer's fixed point theorem.

\proclaim Lemma 1.2. let $Y_0$ be an infinite codimensional
subspace of the Banach space $Y$ and $B$ a compact subset of
$Y$.  Then for every $\t>0$ there is a $y$ in $Y$ with
$\nm{y}<\t$, so that $d(B+y,Y_0) \ge \t/2$.

\proclaim Lemma 1.3. Let $X_0$ be a finite codimensional
subspace of the Banach space $X$.  For every $\t>0$ there is a
compact subset $A$ of $\t \Ball(X)$ so that whenever $\phi$ is
a continuous map from $A$ into $X$ for which
$\nm{\phi(x)-x}< \t/2$ for all $x$ in $A$, then
$\phi(A)\cap X_0 \not= \emptyset$.

\noindent {\bf Proof of Lemma 1.3.\ \ } By the Bartle-Graves
theorem or Michael's selection theorem [Mic], there is a
continuous selection $f: {3\over 4}\t \Ball(X/X_0) \to \t
\Ball(X)$ of the inverse of the quotient mapping
$Q$ from $X$ to $X/X_0$. Set
$A=f\left[{3\over 4}\t \Ball(X/X_0)\right]$, and apply
Brouwer's theorem to the mapping
$x\mapsto x-Q\phi f (x)$ from
$
{3\over 4}\t \Ball(X/X_0)
$ to itself. \hfill\qed

\smsk

\noindent {\bf Proof of Theorem 1.1. \ \ } Get the set $A$ from
Lemma 1.3 for the value $\t = d/2$.  Applying Lemma 1.2 to
the compact set $U[A]$, we get $y$ in $Y$ with $\nm{y} <2b$ so
that
$
d(U[A] + y, Y_0) \ge b
$.
The mapping from $A$ into $X$ defined by \
$
a\mapsto V(Ua+y)
$ \
moves each point of $A$ a distance of at most
$
\omega(V,\nm{y})\le \omega(V,2b)
$; so if this is less than $d/4 = \t/2$, we have from Lemma 1.3
that there must be a point $a_0$ in $A$ for which
$
V(Ua_0+y)
$
is in $X_0$ and hence in $d \Ball(X_0)$.  But then
$Ua_0+y$ would be in $Y_0 + b \Ball(Y)$, which contradicts the
choice of $y$.  Therefore
$\omega(V,2b) \ge d/4$, as desired. \hfill\qed

\smsk

Recall that an unconditional basis is said to
have an upper $p$-estimate (respectively, lower
$p$-estimate) provided that there is a constant $0<C<\infty$ so
that for every finite sequence $\{x_k\}_{k=1}^n$ of vectors which
are disjointly supported with respect to the  unconditional basis,
the quantity $\displaystyle \nm{\sum_{k=1}^n x_k}^p$ is less than
or equal to (respectively, greater than or equal to)
$\displaystyle C^p\sum_{k=1}^n \nm{x_k}^p$.  To characterize the
uniform homeomorphs of $\ell_p$ and the Tsirelson spaces, we need
the following implementation  of the Gorelik Principle.

\proclaim Theorem 1.4.  Suppose that $X$ has an
unconditional basis which has an upper $p$-estimate
and $X$ is uniformly homeomorphic to $Y$.
Then no quotient of $Y$ can have an  unconditional
basis which has a lower $r$-estimate with $r < p$.

For the proof of Theorem 1.4, we need to recall the concept of
{\sl approximate metric midpoint,\/} which plays a role also
in the proofs of Enflo [Enf] and Bourgain [Bou].  Given points
$x$ and $y$ in a Banach space $X$ and $\d\ge 0$, let
$$
\Mid(x,y,\d)=
\left\{ z\in X : \nm{x-z}\vee \nm{z-y}\le (1+\d) {\nm{x-y}\over
2}
\right\}.
$$

We also need a quantitative way of
expressing the well-known fact that a uniformly continuous
mapping from a Banach space is ``Lipschitz for large distances".
If $U$ is a uniformly continuous mapping from a Banach space $X$
into a Banach space $Y$, set for each $t>0$
$$ u\ss{t}= \sup \left\{ {\nm{Ux_1-Ux_2}\over  {t\vee
\nm{x_1-x_2}}} : x_1, \;  x_2 \in X \right\}.
$$

The statement that $U$ is {\sl Lipschitz for large distances\/} is
just that $u\ss{t}$ is finite for each $t>0$.  Obviously
$u\ss{t}$ is a decreasing function of $t$; denote its limit by
$u_\infty$.

\ms
\proclaim Lemma 1.5. Let $U$ be a uniformly continuous mapping
from a Banach space $X$ and let $d>0$.  Suppose that $x$, $y$ in
$X$ satisfy for certain $\e\ge 0$, $\d\ge 0$
\item{(i)} $\nm{x-y}\ge {{2d}\over {1+\d}}$\vskip0.2pt
\item{(ii)} $\nm{Ux-Uy} \ge {u\ss{d} \over (1+\e)} \nm{x-y}$.
\vskip0.2pt\noindent
Then $U\left[\Mid(x,y,\d)\right]\subset \Mid(Ux, Uy,\e+\d+\e\d)$.

\Proof  Let $z$ be in $\Mid(x,y,\d)$.  Then
$${\eqalign{
\nm{Ux-Uz}\vee \nm{Uz-Uy} & \le u\ss{d} \left(d\vee (1+\d)
{\nm{x-y}\over 2}\right) = u\ss{d}  (1+\d){\nm{x-y}\over
2}
\cr & \le
(1+\e)(1+\d){\nm{Ux-Uy}\over2}.\quad\quad\quad\quad\quad
\quad\quad\quad\quad\quad\quad\quad\quad\qed}}
$$

\bs
In a space  which has an unconditional basis which has an upper
or lower
$p$-estimate, the set of
approximate  metric midpoints between two points has some obvious
structure.  We state the next lemma only for points symmetric
around $0$; by translation one obtains a similar statement for
general sets
$\Mid(x,y,\d)$.

\proclaim Lemma 1.6. Suppose that $x$ is in $X$ and $X$ has a
basis $\xn$ whose unconditional constant is one.
\item{(i)} If $\xn$ has an upper $p$-estimate with constant one,
then for each $\d>0$ there is a finite codimensional subspace
$X_0$ of $X$ so that
$\Mid(x,-x,\d\nm{x})\supset \d^{1\over p}\nm{x}\Ball(X_0)$.
\item{(ii)} If $\xn$ has a lower $r$-estimate with constant one,
then for each $\d>0$ there is a finite dimensional subspace $X_1$
of $X$ so that
$\Mid(x,-x,\d\nm{x})\subset X_1+ (r\d)^{1\over r}\nm{x} \Ball(X)
$.

\Proof We assume $p>1$ and $x\not= 0$ since otherwise the
conclusion is trivial. To prove (i), suppose that $y$ is in $X$,
$\nm{y}\le \d^{1\over p}\nm{x}$, and $y$ is disjoint from $x$
relative to the basis $\xn$.  Then
$\nm{y\pm x}\le (1+\d)^{1\over p}\nm{x}< (1+\d)\nm{x}$.
So when $x$ is finitely supported relative to $\xn$, we can
take for  $X_0$ the subspace of all vectors in $X$ which vanish
on the support of $x$.  The general case follows by approximating
$x$ by a vector with finite support; the degree of approximation
depending on $\d$.

The proof of (ii) is similar. If $x$  is  in
$\spa\{x_k\}_{k=1}^n$, then $X_1= \spa\{x_k\}_{k=1}^n$
``works" independently of $\d$; the general case follows by
approximation.\hfill\qed

\ms

{\bf Proof of Theorem 1.4 \ } If the conclusion is false, we may
assume, after renorming  $X$ and
$Y$,  that the  basis $\xn$ for $X$
has unconditional constant one and
also  the upper
$p$-estimate constant for $\xn$ is one, and that there is a
quotient mapping
$Q$ from $Y$ onto a space $Z$ having a  basis $\zn$
with unconditional constant one and with  lower
$r$-estimate  constant one.
Let $U$ be a uniform
homeomorphism from $X$ onto $Y$. Since it is easy to check that
the ``Lipschitz for large distance" constants
$s_t$ of $S=QU$ are bounded away from $0$, we can assume without
loss of generality that $s_\infty=1$.

Fix $\d>0$; later we shall see how small $\d$ need be to yield a
contradiction.  Since $s_t\downarrow 1$ as $t\uparrow \infty$, we
can find a pair $x$, $y$ of vectors in $X$ with $\nm{x-y}$ as
large as we please (for one thing we want $\d^{1\over r}\nm{x-y}
>2$), so that ${\nm{Sx-Sy} \over \nm{x-y}}$ is as close to one as
we please.  From Lemma 1.5 we then get, as long as $\nm{x-y}$ is
sufficiently large, that
$$
S\left[\Mid(x,y,\d \nm{x-y})\right] \subset \Mid(Sx,Sy,2 \d
\nm{x-y}).
\leqno{(1.1)}
$$

By making translations, we only need to consider the case when
$y=-x$ and $U(-x)=-Ux$; that is , we can assume that there is $x$
in $X$ with $\nm{x}$ as large as we like and
$$
S\left[\Mid(x,-x,2\d \nm{x})\right] \subset \Mid(Sx,-Sx, 4 \d
\nm{x-y}).
\leqno{(1.2)}
$$

From Lemma 1.6 we get a finite codimensional subspace $X_0$ of
$X$ and a finite dimensional subspace $Z_0$ of $Z$ so that
$$
S\left[(2\d)^{1\over p} \nm{x} \Ball(X_0)\right] \subset
Z_0 + (4r\d)^{1\over r} \nm{x} \Ball(Z).\leqno{(1.3)}
$$

Set $Y_0=Q^{-1} Z_0$.  Then $Y_0$ has infinite codimension in $Y$
and
$$
U\left[(2\d)^{1\over p} \nm{x} \Ball(X_0)\right] \subset
Y_0 + (4r\d)^{1\over r} \nm{x} \Ball(Y).\leqno{(1.4)}
$$

The Gorelik Principle tells us that
$$
\omega(V,2(4r\d)^{1\over r} \nm{x}) \ge
{{(2\d)^{1\over p} \nm{x}} \over 4},\leqno{(1.5)}
$$
 where $V=U^{-1}$.

  However,
keeping in mind that $\d^{1\over r}2\nm{x } >2$, we have
$$
\omega(V,2(4r\d)^{1\over r} \nm{x}) \le
 v_1 2 (4r\d)^{1\over r} \nm{x}.\leqno{(1.6)}
$$

Putting (1.5) and (1.6) together gives
$$
{{(2\d)^{1\over p}} \over 4} \le
v_1 2 (4r\d)^{1\over r},
$$
which is a contradiction for small enough $\d$. \hfill\qed

\ms

In Section 2 we need the next corollary.

\proclaim Corollary 1.7. Suppose that $X$ has an
unconditional basis and $X$ is uniformly homeomorphic to
$Y$.  Assume that either the unconditional basis for $X$
 has a lower
$r$-estimate for some $r<2$, or that $X$ is superreflexive
and the unconditional basis for $X$ has an upper
$p$-estimate for some $p>2$. Then $\ell_2$ is not isomorphic to
a subspace of $Y$.

\Proof Bourgain [Bou] proved that for $r<2$ there is no
homeomorphism $U$ from $\ell_2$  into $\ell_r$ for which both $U$
and $U^{-1}$ are ``Lipschitz for large distances"; a nonessential
modification yields the same result for any space $X$ which
has an  unconditional basis which  has a lower $r$-estimate for
some $r<2$.  Consequently, $\ell_2$ is not isomorphic to a
subspace of any space which is uniformly homeomorphic to such an
$X$.  So
assume that $X$ is superreflexive and
the unconditional basis for $X$ has an upper
$p$-estimate for some $p>2$.  This implies that $X$ is of type
$2$ in the sense of Maurey-Pisier (see [LT2,~1.f]), and
hence so is $Y$, since
$Y$ is finitely crudely representable in $X$ by Ribe's theorem
[Rib1] (or see [Ben2]).
Thus by Maurey's theorem [Mau3], any
isomorphic copy of $\ell_2$ in $Y$ is complemented in $Y$, so
Theorem 1.4 implies that $\ell_2$ is not isomorphic to a
subspace of $Y$.  \hfill\qed

\bs

\noindent{\bf 2. Uniform Homeomorphs of $\ell_p$, $\Tp$, and
related spaces}

\bs

The first result in  this section is an immediate
consequence of Corollary 1.7 and previously known  results.

\proclaim Theorem {2.1}.  If $X$ is a Banach space which is
uniformly homeomorphic to $\ell_p$, $1<p<\infty$, then $X$ is
isomorphic to
$\ell_p$.

\Proof  By Ribe's theorem [Rib2], $X$ is
isomorphic to a complemented subspace of $L_p$.  But $X$ does not
contain an isomorphic copy of $\ell_2$ by Corollary 1.7 if $p\not =
2$, so by the results of [JO], $X$ is isomorphic to
$\ell_p$. \hfill\qed

\ms

\noindent{\bf Remark.\/}  The case $p=2$ in Theorem 2.1 was proved
by Enflo [Enf].  A simple proof using ultraproducts was provided
by Heinrich and Mankiewicz [HM], [Ben2].  They also used
ultraproducts to give a simple proof of Ribe's theorem [Rib2].
Since we need the type  of reasoning introduced in [HM] to prove the
remaining theorems, we recall some more-or-less standard facts about
ultraproducts of Banach spaces. Much more information, as well as
references to the original sources, can be found in [Hei].

\ms
Given a sequence $\Xn$ of Banach spaces and a free ultrafilter
$\U$ on the natural numbers, denote by $\UXn$ (or just $\UX$ if
all the $X_n$'s are the same) the Banach space ultraproduct of
$\Xn$, defined as the collection of all bounded sequences $\xn$ with
$x_n\in X_n$ under the norm
$
\dstyle \lim_{\U} \nm{x_n}
$.
Here we identify two sequences $\xn$ and $\yn$ as being
the same if
$
\dstyle \lim_{\U} \nm{x_n-y_n}=0
$.
The space $\UX$ is finitely representable in $X$.
A less commonly used fact which is  important for
us is that finite dimensional complemented subspaces of $\UX$
pull down to $X$. Precisely, if $E$ is a finite dimensional
subspaces of $\UX$ which is $\lam$-complemented in $\UX$, then for
each $\e>0$, $E$ is $(1+\e)$-isomorphic to a subspace of $X$ which is
$(\lam+\e)$-complemented in $X$.  This follows from the
{\sl ultrapower version of local duality\/} ([Hei],~p.~90]), which
says that every finite dimensional subspace of
$(\UX)^*$
$(1+\e)$-embeds into
$(X^*)\ss{\U}$
in such a way that the action on any fixed finite subset of
$\UX$ is preserved.  (The space $(X^*)\ss{\U}$ is always a subspace
of $(\UX)^*$, but these spaces coincide only when $X$ is
superreflexive.)

The space $X$ is naturally embedded into $\UX$ as the diagonal.  If
$X$ is reflexive, it is norm one complemented in $\UX$ (map
$\xn$ in $\UX$ to the ${\U}$-weak limit in $X$ of $\xn$); the
kernel $\UXz$ of this projection consists of those bounded
sequences in $X$ which tend weakly to zero  along the ultrafilter
$\U$.

An ultrapower of a Banach lattice $X$ is again a Banach lattice;
moreover, if $X$ is a  $L_p$-space, then the norm in $X$ is
$p$-additive for disjoint vectors, and then so is the norm in the
ultrapower.  Thus (by the generalized Kakutani representation
theorem) the ultrapower is also a
$L_p$-space if $1\le p<\infty$.

Heinrich and Mankiewicz [HM] (see also [Ben2])  used ultraproducts to
give simple proofs of a number of previously known theorems
concerning uniformly homeomorphic Banach spaces and to answer a
number of open problems. Here we just recall the basic approach in
a situation which is general enough to meet our needs.

Suppose that $X$ and $Y$ are separable, uniformly homeomorphic
Banach spaces.  Using the fact that a uniformly continuous mapping
from a Banach space is Lipschitz for large distances, one sees
easily that $X$ and $Y$ have Lipschitz equivalent ultrapowers
$\UX$ and $\UY$.  Suppose now that $X$, hence also $\UX$, is
superreflexive.  Using a back-and-forth procedure and a
classical weak compactness argument [Lin2], one gets a separable,
norm one complemented subspace
$X_1$ of $\UX$ which contains $X$ and is Lipschitz equivalent to a
subspace $Y_1$ of
$\UY$ which contains $Y$.  A differentiation argument (combined
with a technique from [Lin1]) now yields that $Y_1$ isomorphically
embeds into
$X_1$ as a complemented subspace.  This shows that $Y$ is also
superreflexive, and so $X_1$  embeds into $Y_1$ as a complemented
subspace.  Moreover, now that we know that $\UY$ is reflexive,
we can make the earlier construction produce
a $Y_1$ which is norm one complemented in $\UY$.

One of
several consequences of this construction which we use later is
that uniformly homeomorphic spaces $X$ and $Y$ have the same finite
dimensional complemented subspaces if $X$ is separable and
superreflexive (in fact, without any restriction on $X$, but this
requires more work).

\ms

Given the necessary background on ultrapowers,
analysis only slightly more involved than that of Theorem
{2.1} yields:

\proclaim Theorem {2.2}. Let $1<p\ss{1}<p\ss{2}<\dots<p\ss{n} <2$
or
$2<p\ss{1}<p\ss{2}<\dots<p\ss{n}  <\infty$ and set
$\displaystyle X=\sum_{k=1}^n \lpk$.  If $X$ is uniformly
homeomorphic to a Banach space $Y$, then $X$ is isomorphic to $Y$.

\Proof  For simplicity of notation, we treat the case $n=2$.
By one of the results of
Heinrich-Mankiewicz [HM] (or see [Ben2]), $X$, $Y$ have Lipschitz
equivalent ultrapowers $\UX$,
$\UY$, respectively; and $\UX$ splits as the direct sum of an $\Lpo$
space with an
$\Lpt$ space:
$\UX = \Lpo(\tilde\mu)\oplus  \Lpt(\tilde\mu)$. (We do not need that
the measures for $p_1$ and $p_2$ are the
same, but they are.)  Ribe proved that $Y$ is
finitely crudely representable in $X$ ([Rib1], [HM], [Ben2]), so $Y$
is superreflexive and
$\UY$ splits as the direct sum of $Y$ and some space $\tilde Z$
(which is, incidentally, embeddable as a subspace of
$
\Lpo(\nu)\oplus  \Lpt(\nu)$ for some measure
$\nu$).  Since $Y$ is separable, there exists a separable subspace
$Z$ of  $\tilde Z$ so that the image of $Y+ Z$ under the Lipschitz
isomorphism from $\UY$ onto $\UX$ is of the form
$\Lpo(\mu)\oplus \Lpt(\mu)$, where $\mu$ is the restriction of
$\tilde\mu$ to a separable sigma subalgebra.  Hence by [HM], $Y$ is
isomorphic to a complemented subspace of $\Lpo(\mu)\oplus
\Lpt(\mu)$, whence of $\Lpo(0,1)\oplus \Lpt(0,1)$.

Suppose now that $\po$ and $\pt$ are larger than two.  Let $J$
denote the embedding of $Y$ onto a complemented subspace of
$\Lpo(0,1)\oplus \Lpt(0,1)$ and let $S_1$, $S_2$ be the natural
projections from $\Lpo(0,1)\oplus \Lpt(0,1)$ onto $\Lpo(0,1)$ and
$\Lpt(0,1)$, respectively.  We know from Corollary 1.7 that $\ell_2$
does not embed into $Y$, so by the generalization in [Joh1] of the
theorem from [JO] used earlier, for $k=1$, $2$, the operator
$S_kJ$ factors through $\lpk$.  Hence  $J$ factors through
$\lpo\oplus\lpt$, whence $Y$ embeds into $\lpo\oplus\lpt$ as a
complemented subspace.  However, by a result of
Edelstein-Wojtaszczyk [Ede], [Woj], [EW] (or see [LT1,~p.~82]), every
complemented subspace of
$\lpo\oplus\lpt$ is isomorphic to $\lpo$, $\lpt$, or
$\lpo\oplus\lpt$.   Theorem 2.1 eliminates the first two
possibilities.

When $\po$ and $\pt$ are smaller than two, we pass to the duals:
$Y^*$ is isomorphic to a complemented subspace of $\Lqo(0,1)\oplus
\Lqt(0,1)$, where $\qk$ is the conjugate index to $\pk$.   $Y^*$
also does not contain an isomorph of $\ell_2$ (since  $Y^*$ has
type $2$, every copy of $\ell_2$ in $Y^*$ is complemented).  The
reasoning in the last paragraph shows that $Y^*$ is isomorphic to
$\lqo\oplus\lqt$.\hfill\qed

\ms

In sections 4 and 5 we prove that for $1\le p <\infty$, the
$p$-convexified version $\Tp$ of Tsirelson's space is uniformly
homeomorphic to $\Tp\oplus \ell_p$.  Here we show that for $ 1<
p<\infty$,  $\Tp$  is uniformly homeomorphic to at most two
isomorphically distinct spaces.  Since it is desirable to have
general conditions which limit the isomorphism class
of spaces uniformly homeomorphic to a given space, it seems
worthwhile to  prove some
results  in this direction which might be used elsewhere. The reader
who is interested only in the example can skip to
Proposition 2.7 and substitute $\Tpk$ for $X_k$ in the
statement.

We begin with a definition.
\ms
\noindent {\bf Definition {2.3.}} A
Banach space $X$ is said to be {\sl
as.\/} ${\cal L}_p  $, $1 \le p \le
\infty$, provided there exists $\lam$ so
that for every $n$ there is a finite
codimensional subspace $Y$ so that every
$n$ dimensional subspace of $Y$ is
contained in a subspace of $X$ which is
$\lam$-isomorphic to $L_p(\mu)$ for some $\mu$.
\ms
A space is as. $\sLt  $ if and only if it is asymptotically
Hilbertian in the sense of Pisier [Pis].   We avoided the term
asymptotically ${\cal L}_p  $ because  when
$p\not=2$ as.
${\cal L}_p  $ may not be  the
``right" definition for asymptotically
${\cal L}_p$.   The definition we give is also not the definition one
gets by specializing Pisier's [Pis, p. 221] {\sl as. property\/}
$(P)$ to
$P = {\cal L}_p$-structure; however, Pisier's definition seems
``right" only for hereditary properties.

In section 5 we review some of the properties of the Tsirelson
spaces.  Here we just mention that
$\Tp$ has an unconditional basis
$\en$ for which there is a constant $\lam$ such that for every $n$,
there exists an $m$, so that every $n$-tuple of disjointly supported
unit vectors in $\spa\{e_k\}_{k=m}^\infty$ is $\lam$-equivalent to
the unit vector basis in $\ell_p^n$.  It is easily seen that a space
$X$ with this property is  as. ${\cal L}_p $ and
also that $\UXz$ is isomorphic to $L_p(\mu)$ for some measure $\mu$.
This property of such spaces can be generalized to the class of as.
${\cal L}_p$ spaces.

\proclaim Proposition {2.4.a}.  Suppose that $X$ is as. ${\cal L}_p$,
$1<p<\infty$.  Then for every ultrafilter $\U$ on the natural
numbers, the space $ \UXz$ is a ${\cal L}_p$ space.

\Proof  First use the James' characterization of
reflexivity/superreflexivity to see that
$X$ is superreflexive (see, for example, the argument in
[Pis,~pp.~220\&222] for the similar case of asymptotically
Hilbertian spaces).  Therefore
$\UX$ splits as the direct sum of $X$ and $\UXz$.  Now if $\tilde E$
is a finite; say, $n$; dimensional subspace of $\UXz$, then it is not
hard to see that there exists a sequence $\Ek$ of $n$-dimensional
subspaces of
$X$ which converges weakly to zero along $\U$ (in the sense that every
bounded sequence $\xk$ with $x_k$ in $E_k$ converges weakly to $0$
along $\U$) so that
$\UEk = \tilde E$.  Since $\Ek$ converges weakly to zero along $\U$ and
have bounded dimension, given any finite codimensional subspace $Y$ of
$X$, a standard small perturbation argument shows that there are space
automorphisms $T_k$ on $X$ with
$\displaystyle \lim_{\U}\nm{I-T_k} = 0$ and
$T_k E_k \subset Y$ for all $k$.  Consequently, since $X$
is as. ${\cal L}_p$ (with constant $\lam$, say), there are superspaces
$F_k$ in $X$  of
$E_k$ with $F_k$ $\lam_k$-isomorphic to $L_p(\mu_k)$ and
$\displaystyle \lim_{\U} \lam_k \le \lam$.  In particular,
the inclusion mappings from the $E_k$'s into $X$ uniformly factor
through $L_p$-spaces for a $\U$-large set of $k$'s.  Passing to
ultraproducts and using the fact that an ultraproduct of $L_p$-spaces
is again an
$L_p$-space, we conclude that the inclusion mapping from $\tilde E$
into $\UX$
 factors through an $L_p$-space with the
norms of the factoring maps independent of  $k$.  The same
can be said about the inclusion mapping from  $\tilde E$ into $\UXz$
since $\UXz$ is complemented in $\UX$.  Finally, we conclude from
[LR] that $\UXz$ is either a ${\cal L}_p$ or isomorphic to a Hilbert
space.  But the latter is impossible for $p\not= 2$. Indeed,
$\ell_p$ is  finitely crudely representable in $X$ and hence also in
every finite codimensional subspace of $X$, and this implies that
$\ell_p$ embeds into $\UXz$. \hfill\qed
\ms
\proclaim Proposition {2.4.b}.  Let $X$ be a
reflexive space for which there exists  an ultrafilter $\U$ on the
natural numbers so that $\UXz$ is $\sLp$, $1 < p <\infty$. Then $X$ is
as. $\sLp$.

\Proof Assume that the $\sLp$ constant of $\UXz$ is smaller than
$\lam$, but that $\UXz$ is not as. $\sLp$ with constant
$\lam$. Let $n$ be such that every finite codimensional subspace
of $X$ contains an $n$-dimensional subspace no superspace of which
is $\lam$-isomorphic to an $L_p$ space.  Using a standard basic
sequence  construction, we can construct a finite dimensional
decomposition
$\Ek$ of
$n$-dimensional subspaces of $X$ so that no superspace of any $E_k$
is $\lam$-isomorphic to an $L_p$ space.  But since $X$ is reflexive,
the sequence $\Ek$ converges weakly to zero, and hence the
ultraproduct $E=\UEk$ is an $n$-dimensional subspace of $\UXz$.
Then there is a finite (say, $m$,) dimensional subspace
of $\UXz$ which contains $E$ and whose Banach-Mazur
distance from $\ell_p^m$ is less than $\lam$.  It
follows that $F$ can be represented as
$\UFk$, where for each $k$ $F_k$ is an $m$-dimensional subspace of $X$
which contains $E_k$.  Necessarily the set of $k$'s for which $F_k$
is $\lam$-isomorphic to $\ell_p^m$ is in $\U$, hence is nonempty.
\hfill\qed

\ms

The next lemma is an as. version of the result [LR] that a
complemented subspace of a
$\sLp$ space is either a $\sLp$ space or a $\sLt$ space.

\proclaim Lemma {2.5}. If $X$ is as. $\sLp$, $1<p<\infty$, and $Y$
is a complemented subspace of $X$, then $Y$ is as. $\sLp$ or as.
$\sLt$.  $Y$ is as. $\sLp$ iff $\ell_p$ is finitely crudely
representable in $Y$ iff $Y$ contains uniformly complemented
$\ell_p^n$'s for all $n$.

\Proof The ultraproduct of a projection from $X$ onto $Y$ defines a
projection from $\UX$ onto $\UY$ which projects $\UXz$ onto $\UYz$.
Thus the first conclusion follows from Proposition 2.4  and the
classical theory of $\sLp$-spaces [LPe], [LR].  The last sentence
follows from the fact [KP] that a complemented subspace of an
$L_p$ space which is not isomorphic to a Hilbert space contains a
complemented subspace isomorphic to $\ell_p$ and the fact,
mentioned earlier, that one can pull down finite dimensional
well-complemented subspaces of an ultrapower to the base space.
\hfill\qed

\ms

We also need the as. version of Ribe's [Rib2] theorem that a uniform
homeomorph of  a
$\sLp$ space,  $1<p<\infty$, is again a $\sLp$ space.

\proclaim Lemma {2.6}. If $X$ is  separable as. $\sLp$,
$1<p<\infty$, and $Y$ is uniformly homeomorphic to $X$, then $Y$ is
as. $\sLp$.

\Proof  $Y$ is finitely crudely representable in $X$ by [Rib1]
(or see [HM]), so
$Y$ is superreflexive. Reasoning as in the proof of Theorem
2.2, we get Lipschitz equivalent ultraproducts $\UX=X\oplus \UXz$
 and $\UY=Y\oplus \UYz$ of $X$ and $Y$, respectively; and, by
Proposition 2.4, $\UXz$ is $\sLp$.  Continuing as in  the earlier
proof, we get separable subspaces $X_1$ of $\UXz$ and $Y_1$ of
$\UYz$ with $X_1$ a $\sLp$ space so that $X\oplus X_1$ is Lipschitz
equivalent to  $Y\oplus Y_1$.  By another result from [HM], $Y$ is
isomorphic to a complemented subspace of $X\oplus X_1$, which of
course is as. $\sLp$.  Hence by Lemma 2.5, $Y$ is as. $\sLp$ or
as. $\sLt$.  But $X$ is finitely representable in $Y$, so the
latter is impossible.\hfill\qed

\ms

\noindent{\bf Remark.  \  }  The separability assumption is not
needed in Lemma 2.6.  One way to see this is to check
that if
$X$ is as. $\sLp$, $1<p<\infty$, then $X=X_1\oplus X_2$ with $X_1$
separable and $X_2$ $\sLp$.  Casazza and Shura [CS,~p.~150] proved
this for
$p=2$ and the general case can be done similarly.

\ms

The proof of Lemma 2.6 yields that if the separable as. $\sLt$-space
$X$ is uniformly homeomorphic to $Y$, then $X\oplus\ell_2$ and
$Y\oplus\ell_2$ are Lipschitz equivalent and hence by [HM] each is
isomorphic to a complemented subspace  of the other.  Proposition
2.7 gives a version of this for certain  direct sums of as. $\sLp$
spaces. For its proof, we need a  consequence of the
preceding. Suppose that $X\oplus Y$ contains uniformly
complemented copies of
$\ell_r^n$ for all $n$,  $X$ is as. $\sLp$,  $1<p<\infty$, and
$2\not= r \not= p$.  Then $Y$ contains uniformly complemented copies
of $\ell_r^n$ for all $n$.  Indeed, one observes that $\ell_r$
embeds complementably into
$\left(X\oplus Y\right)\ss{\U}=\UX\oplus\UY$, and hence
complementably into either $\UX$ or $\UY$.  But $\UX$ is as.
$\sLp$ (for example, by Proposition 2.4.a), so if $\ell_r$ embeds
complementably into $\UX$, then $\ell_r$ is as. $\sLp$ or as.
$\sLt$ by Lemma 2.5.  The latter is clearly impossible, and so is
the former by the last statement in Lemma 2.5 and well-known
properties of $\ell_r$.  Hence $\ell_r$ embeds complementably into
$\UY$, whence $\UY$, {\sl a fortiori\/} $Y$ itself, contains
uniformly complemented copies of $\ell_r^n$ for all $n$.

\proclaim Proposition {2.7}. Let
$1<p\ss{1}<p\ss{2}<\dots<p\ss{n} <2$ or
$2<p\ss{1}<p\ss{2}<\dots<p\ss{n}  <\infty$ and for $1\le
k\le n$, assume that $X_k$ is a separable, as.
$\sLpk$ space with an unconditional basis which has an
lower $r$-estimate for some $r<2$ or an upper
$p$-estimate for some $p>2$.
 Set
$\displaystyle X=\sum_{k=1}^n X_k$.  If
$X$ is uniformly
homeomorphic to a Banach space $Y$, then $X\oplus \sum_{k=1}^n \lpk$
and
$Y\oplus  \sum_{k=1}^n \lpk$ are  isomorphic to complemented
subspaces of each other.

\Proof  As in the proof of Theorem 2.2, we assume, for simplicity
of notation, that $n=2$.  Reasoning as in the proof of Lemma
2.6, we get Lipschitz equivalent ultraproducts $\UX$ and $\UY$
of $X$ and $Y$, respectively.  $\UX$ is
$X\oplus \tilde {X_1} \oplus \tilde {X_2}$, where
$\tilde {X_k}$ is $\sLpk$ for $k=1,2$. Continuing as in the proof of
Lemma 2.6, we conclude that $X\oplus Z_1 \oplus Z_2$ is Lipschitz
equivalent to $Y \oplus W$ for some separable $\sLpk$
spaces $Z_k$ and some  $W$.  By enlarging the $Z_k$'s and
$W$, if necessary, we can assume that $Z_k=\Lpk(0,1)$.  So $Y$
embeds into
$\left(X_1\oplus \Lpo\right)\oplus \left(X_2\oplus \Lpt\right)$ as a
complemented subspace.  However, we also know from Corollary 1.7 that
$\ell_2$ is not isomorphic to a subspace of $Y$, so, reasoning as in
the proof of Theorem 2.2, we conclude that $Y$ is isomorphic to a
complemented subspace of
$\left(X_1\oplus \lpo\right)\oplus \left(X_2\oplus \lpt\right)$.
The two spaces in parentheses are totally incomparable, so by the
theorem of Edelstein-Wojtaszczyk [EW] (or [LT1,~p.~80]), $Y=Y_1\oplus
Y_2$ with
$Y_k$   isomorphic to a complemented subspace of
$X_k\oplus \lpk$.  Lemma 2.5  now tells us that $Y_k$ is as.
$\sLpk$ or as. $\sLt$.  But $Y$ is uniformly homeomorphic to $X$ and
$X$ contains uniformly complemented $\lpk^n$'s for all $n$,
hence so does $Y$ by [HM] (see the earlier discussion on
ultraproducts). But
$Y_1$ cannot contain  uniformly complemented $\lpt^n$'s for all $n$,
 and so $Y_2$ must.  Thus $Y_2$ is as.
$\sLpt$ and, similarly, $Y_1$ is as.  $\sLpo$. This shows that  the
situation with $X$ and
$Y$ is  symmetrical, and we can conclude that $X$ is isomorphic
to a complemented subspace of $Y\oplus  \lpo \oplus  \lpt$.
\hfill\qed
\ms
\noindent{\bf Remark.\/} In Proposition 2.7, even when $n=1$ the
hypothesis on
$X_1$ cannot be weakened to ``$X$ is
$\sLpo$,
 $X$ has an unconditional basis, and no subspace
of $X$ is isomorphic to $\ell_2$".  Indeed, $\Tpo\oplus \Tt$ is
as.  $\sLpo$ for $1<p_1<\infty$ and thus is a counterexample by
Proposition 5.7.

\ms
The conclusion in Proposition 2.7 is not completely satisfactory,
since it is open whether a complemented subspace of
$\Tp$ which contains a complemented copy of $\Tp$ must be
isomorphic to $\Tp$.  Fortunately, it is possible to improve
Proposition 2.7 for spaces of Tsirelson type without solving this
problem about $\Tp$.  However, this requires rather more work,
so we treat first the simpler case of $\Tt$ and related spaces,
which already provides examples of spaces that  are isomorphic
to exactly two isomorphically distinct spaces.  Basic to this
argument is a recent result concerning differentiation [LP]:

\proclaim \phantom. Assume that $X$ and $Y$ are separable
superreflexive Banach spaces.  Let $F$ be a Lipschitz map from $X$
into $Y$ and $G$ a Lipschitz map from $X$ into a finite
dimensional normed space.  Then for every $\e>0$, there is a point
$x_0$ in $X$ so that $F$ and $G$ are Gateaux differentiable at
$x_0$ and moreover there exists $\d>0$ so that for $||u||\le\d$,
$$
||G(x_0+u)-G(x_0)-S u||\le \e ||u||,
$$
where   $S$ is the Gateaux derivative of $G$ at $x_0$.

In other words, $G$ is ``almost" (up to $\e$) even Frech\'et
differentiable at $x_0$.  An easy consequence of this result and
the reasoning used in [HM] is the following:

\proclaim Proposition {2.8}. Suppose that $F$ is a mapping from the
separable, superreflexive space $X$ onto a Banach space
$Y$, the Lipschitz constant of $F$ is one and the Lipschitz constant
of $F^{-1}$ is $C<\infty$ .  Then  for every finite dimensional
subspace $E$ of
$Y^*$ and each $1>\e'>0$, , there is a norm one operator $T$ from $X$
into
$Y$ so that
$\nm{T^{-1}}\le C$, $TX$ is $C$-complemented in $Y$, and \
$\nm{{T^*_{|E}}^{-1}}\le (1-\e')^{-1}C$.

\Proof
 Note that the space
 $Y$ is also superreflexive (e.g., by [Rib1]).  We shall see that
for
$T$ one may take the Gateaux derivative of $F$ at a suitable point
$x_0$ in
$X$; the point depends on the subspace $E$ of
$Y^*$.   In view of (the proof of)
Proposition 2.1 of [HM] and the fact that $T$ is the Gateaux
derivative of $F$ at some point, the conclusions that
$\nm{T^{-1}}\le C$ and $TX$ is
$C$-complemented in $Y$ are automatic.  To find the point $x_0$,
apply the result from [LP] to the mapping $F$ and an auxiliary
mapping $G$, which we now define:  Let $G$ be the map $F$
followed by the evaluation (quotient) mapping $Q$ from $Y$ onto
$E^*$, so that  for $x$ in $X$ and $y^*$ in $E$,
$Gx(y^*)=  y^*(Fx)$.  Given $\e'>0$, set $\e={{\e'}\over C}$ and get
$x_0$, $\d>0$,  and
$S$ from the [LP] result stated above. For notational simplicity,
assume that $x_0=0$ and $F(0)=0$  (this  is ``without loss
of generality" because it amounts to making one translation in $X$
and another in
$Y$). Of course, $S=QT$, where $T$ is the Gateaux derivative of $F$ at
$0$.  Only the last conclusion, that
$\nm{{T^*_{|E}}^{-1}}\le (1-\e')^{-1}C$, needs to be checked. Since $Q^*$
is the inclusion mapping from $E$ into $Y^*$, this is the same as checking
that $\nm{S^{*-1}}\le (1-\e')^{-1}C$. So suppose that
$y^*$ is a norm one vector in
$E$ and choose a vector
$y$ in
$Y$ of norm $\d\over C$ so that $y^*(y)$ is equal to $\d\over C$
(an approximation would be OK but is not needed because $Y$ is
reflexive).  Set $x=F^{-1}y$, so that $||x||\le\d$. Then
$$\eqalign {
|S^*y^*(x)-{\d\over C}|&=|S^*y^*(x)-y^*(Fx)|=|y^*(Sx)-Gx(y^*)|
\cr
&\le
||y^*||\, ||Sx-Gx||
\le {\e'\over C} ||x||\le {\e'\over C} \d .\cr}
$$
Since $||x||\le\e$, this gives $||S^*y^*||\ge {{1-\e'}\over C}$, so
that
$\nm{{S^*}^{-1}}\le (1-\e')^{-1}C$, as desired.
\hfill\qed

\ms

\proclaim Proposition {2.9}.  Suppose that $X$ is separable, as.
$\sLt$,
 does not contain a complemented subspace
isomorphic to
$\ell_2$, and is isomorphic to its hyperplanes.  If $X$ is uniformly
homeomorphic to $Y$, then $Y$ is isomorphic either to $X$ or to
$X\oplus \ell_2$.

\Proof The main part of the proof is devoted to proving the following

\proclaim Claim. $X \oplus \ell_2$ is isomorphic to  $Y \oplus \ell_2$.

Assuming the Claim, we complete the proof of Proposition 2.9 as
follows: If $Y$ contains a complemented subspace isomorphic to
$\ell_2$, then $Y$ is isomorphic to $Y \oplus \ell_2$ and so we are
done, so assume otherwise. By a slight abuse of notation, write
$X \oplus \ell_2 =  Y \oplus W$ with $W$ isomorphic to $\ell_2$. Since
$X$ does not contain a complemented subspace isomorphic to $\ell_2$,
every operator from $X$ into
$\ell_2$ is strictly singular; that is, not an isomorphism on any
infinite dimensional subspace of $Y$.  Hence by the result of
Edelstein-Wojtaszczyk [EW], by applying a space automorphism to $X
\oplus \ell_2$, we can assume, without loss of generality, that
$Y=Y\ss{X}\oplus Y_2$, $W=W\ss{X}\oplus W_2$ with $Y\ss{X}$,
$W\ss{X}$ subspaces of $X$ and $Y_2$, $W_2$ subspaces of $\ell_2$.
But $Y_2$ must be finite dimensional since $Y$ does not contain a
complemented subspace isomorphic to $\ell_2$.
Similarly, $W\ss{X}$ is finite dimensional since since $X$ does not
contain a complemented subspace isomorphic to $\ell_2$.  Thus $Y$ is
isomorphic to $X$ plus or minus a finite dimensional space,
which is isomorphic to $X$.

We turn to the proof of the Claim.

As mentioned earlier, the proof of Lemma 2.6 shows that
$Y'\equiv Y\oplus
\ell_2$ is Lipschitz equivalent to $X'\equiv X\oplus \ell_2$.
Let $\En$ be an increasing sequence of finite dimensional subspaces
of $X^*$ whose union is dense in $X^*$ and, using Proposition 2.8,
get for each $n$ a norm one isomorphism $T_n$ from $Y'$ onto a
$C$-complemented subspace of $X'$ with $\nm{T_n^{-1}}\le
C$, $\nm{{T_n^*}_{|E_n}^{-1}}\le C$.  The ultraproduct ${T}$ of
the isomorphisms $T_n$ is an isomorphism from the ultrapower
$\UYp$ onto a complemented subspace of the ultrapower $\UXp$:
${T}(y_1,y_2,\dots)=(T_1y_1,T_2y_2,\dots)$.  Since $X'$ is
superreflexive, $(\UXp)^* = X'^*\ss{\U}$ and  the adjoint
$T^*$ is defined for $\xn$ a bounded sequence in $X'^*$ by
$T^*(x^*_1,x^*_2,\dots)=(T_1^*x^*_1,T_2^*x^*_2,\dots)$.
Identify $X^*$ with
$\{(x^*,x^*,\dots) : x^* \in X^* \}$;
evidently  $T^*_{|X^*}$ is an isomorphism.  By
Proposition 2.4.a,
$X'^*\ss{\U}$ is isomorphically just
$X^*$ direct summed with a nonseparable Hilbert space.  Thus by coming
down to appropriate separable subspaces of the ultrapowers, we get
that there is a projection
$P$ from $X\oplus \ell_2$ onto a subspace isomorphic to
$Y\oplus \ell_2$
for which $P^*_{|X^*}$ is an isomorphism. This implies that the
projection from  $X^* \oplus \ell_2$ onto
$\ell_2$ is an isomorphism on the subspace ${\cal R}(P)^\perp$ of $X^*
\oplus
\ell_2$; in particular, ${\cal R}(P)^\perp$ is isomorphic to a
Hilbert space.  But ${\cal R}(P)^\perp$ is isomorphic to the dual of
the null space of $P$, so the null space of $P$ is also isomorphic to
a Hilbert space.
\hfill\qed

\ms

\noindent{\bf Remark.\/} The last
two hypotheses on $X$ were not used in the proof of the Claim.
\ms
\noindent{\bf Remark.\/} The proof of Proposition 2.9 yields that
if $X$ is a separable as. $\sLt$ space and every infinite
dimensional subspace of
$X$ contains a complemented subspace which is isomorphic to
$\ell_2$, then $X$ is determined by its uniform structure. The
simplest examples of such spaces which are different from $\ell_2$
are $\ell_2$-sums of $\ell_{p_n}^{k_n}$ with $p_n\to 2$ and $k_n\to
\infty$ sufficiently quickly.

\proclaim Proposition {2.10}. Let $1<p\ss{1}<p\ss{2}<\dots<p\ss{n}
<2$ or
$2<p\ss{1}<p\ss{2}<\dots<p\ss{n}  <\infty$ and for
$1\le k\le n$, assume that $X_k$ is a separable, as. $\sLpk$,
has an unconditional basis which has an upper
$p$-estimate for some $p>2$ or a lower $r$-estimate for some
$r<2$, is isomorphic to its hyperplanes, and $X_k^*$ does
not contain a subspace isomorphic to $\ell_s$ for
any $s$. Set
$\displaystyle X=\sum_{k=1}^n X_k$.  If
$X$ is uniformly homeomorphic to a Banach space $Y$, then $Y$ is
isomorphic to $X\oplus
\sum_{k\in F} \lpk$ for some subset $F$ of $\{1,2,\dots,n\}$.

\Proof As usual, assume $n=2$.  From the proof of Proposition 2.7,
we know that $Y=Y_1\oplus Y_2$ with $Y_k$ as. $\sLpk$ and that
$X\oplus \Lpo(0,1)\oplus \Lpt(0,1)$ is Lipschitz equivalent to
$Y\oplus \Lpo(0,1)\oplus \Lpt(0,1)$.  The main part of the
proof is devoted to proving the following

\proclaim Claim. There exists a quotient $W$ of
$\Lpo(0,1)\oplus \Lpt(0,1)$ so that $Y\oplus W$ is isomorphic to
$X \oplus \Lpo(0,1)\oplus \Lpt(0,1)$.

Assuming the Claim, we complete the proof as follows: Every
operator from $X^*$ to  $\Lqo(0,1)\oplus \Lqt(0,1)$  (where
$q\ss{k}$ is the conjugate index to $p\ss{k}$) is strictly
singular by [KP], [Ald], [KM], so the Edelstein-Wojtaszczyk
theorem [EW], [LT1, p. 80] says that by applying a space
automorphism to $X\oplus
\Lpo(0,1)\oplus\Lpt(0,1)$, we may assume, without loss of
generality, that
$Y=Y\ss{X} \oplus Y\ss{L}$, $W=W\ss{X}\oplus W\ss{L}$ with
$Y\ss{X}$, $W\ss{X}$ subspaces of $X$ and $Y\ss{L}$, $W\ss{L}$
subspaces of $\Lpo(0,1)\oplus \Lpt(0,1)$.  As in the proof of
Theorem 2.2, Theorem 1.7 and old results from the linear theory
imply that $Y\ss{L}$ is isomorphic to $\lpo$, $\lpt$,
$\lpo \oplus \lpt$, or is finite dimensional (in which case we
can assume the dimension is zero because $X$ is isomorphic to
its hyperplanes).  So the proof is complete once we
observe that $W\ss{X}$ must be  finite dimensional.  But
$W\ss{X}^* $ embeds into $\Lqo(0,1)\oplus \Lqt(0,1)$
 and hence by
[KP], [Ald], [KM] would contain a subspace isomorphic to
$\ell_s$ for some $s$ if it were infinite dimensional. This
would contradict the hypotheses on $X^*$ since
$W\ss{X}^* $ embeds into $X^*$.

We turn to the proof of the claim, which is only slightly more
complicated than the proof of the Claim in Proposition 2.9.

Since $X'\equiv X\oplus \Lpo(0,1)\oplus \Lpt(0,1)$ is Lipschitz
equivalent to
$Y'\equiv Y\oplus \Lpo(0,1)\oplus \Lpt(0,1)$, we use Lemma 2.6
just as in the proof of Proposition 2.9 to get good isomorphisms
$T_n$ from
$Y'$ onto well-complemented subspaces of $X'$ so that
$T_n^*$ is a  good isomorphism on a   subspace
$E_n$ of $Y^*$, where $\En$ is an increasing sequence of
finite dimensional subspaces of $Y^*$ whose union is dense in
$Y^*$.  Taking ultrapowers, we get an isomorphism $T$ from an
ultrapower ${Y'}\ss{\U}$ of $Y'$ onto a complemented subspace of
${X'}\ss{\U}$ so that
$T^*$ is an isomorphism on $X^*$.  But we know that
${Y'}\ss{\U}$; respectively, ${X'}\ss{\U}$, is just $Y$;
respectively $X$, direct summed with the direct sum of a
$\sLpo$-space with a $\sLpt$-space.  Therefore by coming down to
appropriate separable subspaces of the ultrapowers and, if
necessary, enlarging the $\sLpk$-spaces, we conclude that there
is a projection $P$ from $X\oplus \Lpo(0,1)\oplus \Lpt(0,1)$ onto
a subspace isomorphic to $Y\oplus \Lpo(0,1)\oplus \Lpt(0,1)$ for
which $P^*_{|X^*}$ is an isomorphism.  This implies that the
natural projection from
$X^*\oplus \Lqo(0,1)\oplus\Lqt(0,1)$ onto
$\Lqo(0,1)\oplus\Lqt(0,1)$is an isomorphism on the subspace
${\cal R}(P)^\perp$; in particular, ${\cal R}(P)^\perp$ is
isomorphic to a subspace of $\Lqo(0,1)\oplus\Lqt(0,1)$.  But
${\cal R}(P)^\perp$ is isomorphic to the dual of the null space
of $P$, so the null space of $P$  is isomorphic to a quotient of
$\Lpo(0,1)\oplus \Lpt(0,1)$. \hfill\qed

\bs

\noindent{\bf 3. Uniform Homeomorphs of $c_0$}

\bs
In this section we prove:

\proclaim Theorem {3.1}. Let $X$ be a Banach space which has
$C(\omega^\omega)$ as a quotient space.  Then $X$ is not uniformly
homeomorphic to $c_0$.

In particular, $c_0$ is not uniformly homeomorphic to $C[0,1]$ or any
other $C(K)$ space except those which are isomorphic to $c_0$,
so Theorem 3.1 solves the problem mentioned in [Aha]
whether $c_0$ and $C[0,1]$ are Lipschitz equivalent.
In fact,
an immediate application of Theorem {3.1} and results of
Alspach,  Zippin, and Benyamini   is:

\proclaim Corollary {3.2}. Let $X$ be a complemented subspace of a
$C(K)$ space.  If $ X$ is uniformly homeomorphic to $c_0$, then $X$
is isomorphic to $c_0$.

\Proof  By [Als] (or see [AB]), since $X$ does not have a quotient
isomorphic to $C(\omega^\omega)$, the $\e$-Szlenk index of $X$ is
finite for each $\e>0$.  But in [Ben1, proof of Theorem 3] it was
shown (using, among other things, a small variation on a lemma in
[Zip]) that then $X$ must be isomorphic to a quotient of $c_0$,
hence $X$ is isomorphic to $c_0$ by [JZ].

Although it is known [HM] that a uniform homeomorph of $c_0$ is a
$\sLinf$ space, it is open whether it must be isomorphic to a
complemented subspace of a $C(K)$ space. 

\ms

\noindent{\bf Proof of Theorem {3.1}.\ \ } Denote by $K$ the space
$\omega^\omega$ and let $K^{(1)}$, $K^{(2)}$, \dots be the derived
sets of $K$.  Let $R_n$ be the restriction mapping from $C(K)$ onto
$C(K^{(n)})$.

Assume that there is a uniform homeomorphism $U$ from $c_0$ onto $X$
with inverse $V$ and that there is a quotient mapping $Q$ from $X$ onto
$C(K)$.  Without loss of generality we assume that $Q$ maps the open
unit ball of $X$ onto the open unit ball of $C(K)$.

For each $n$, let $s_{n,\infty }$ be the limiting Lipschitz constant
of the map $S_n\equiv R_nQU$.  So $s_{n,\infty }$ is the
smallest constant so that for each $\e>0$, there is a
$d=d(\e,n)$ so that if $y$, $z$ are in $c_0$ with $\nm{y-z}\ge d$,
then $\nm{S_ny-S_nz}\le (s_{n,\infty }+\e)\nm{y-z}$.  The sequence
$\sninf$ is decreasing and tends to a limit $s_{\infty,\infty}$ which
is easily seen to be positive, so we can assume without loss of
generality that $s_{\infty,\infty}=1$.

Fix $\e>0$ small and take $n_0$ so that $s_{n_0,\infty }< 1+\e$.
We can find a pair $x$, $y$ of vectors in $c_0$ with $\nm{x-y}$ as
large as we please so that
$(1-\e)\nm{x-y}<\nm{S_{n_0+1}x-S_{n_0+1}y}$; and, since $\nm{x-y}$ is
large, of necessity $\nm{S_{n_0}x-S_{n_0}y}< (1+\e)\nm{x-y}$.  By
making translations, we can assume that $y=-x$  and $U(-x)=-Ux$, so
we have:
$$
(1-\e)\nm{x}<\nm{S_{n_0+1}x}\le \nm{S_{n_0}x}< (1+\e)\nm{x}.
\leqno{(3.1)}
$$

Moreover, since  $\nm{x}$ can be taken arbitrarily large, we can
assume, in view of Lemma 1.5 that
$$
S_{n_0}\left[\Mid(x,-x,0)\right] \subset
\Mid(S_{n_0}x,-S_{n_0}x, \e).
\leqno{(3.2)}
$$

Choose $N$ so that the magnitude of the $k$-th coordinate of $x$ is
less than $\e$ for $k>N$ and let $Y_0$ be the finite
codimensional subspace of $c_0$ consisting of those vectors whose
first $N$ coordinates are zero, so  that
$$
\left(\nm{x}-\e\right)\Ball(Y_0)\subset \Mid(x,-x,0).
\leqno{(3.3)}
$$

Now take $t_0$ in $K^{(n_0+1)}$ with
$|S_{n_0+1}x(t_0)|=\nm{S_{n_0+1}x}> (1-\e)\nm{x}$.  The point $t_0$ is
a limit point of $K^{(n_0)}$, so $|S_{n_0}x(t)|> (1-\e)\nm{x}$ for
infinitely many points $t$ in $K^{(n_0)}$; say, for $t\in B$.

Let $W$ be the subspace of $C(K^{(n_0)})$ of functions which vanish
on $B$. Recalling that
$\nm{S_{n_0}x}< (1+\e)\nm{x}$, we observe that
$$
\Mid(S_{n_0}x,-S_{n_0}x, \e)\subset
W+2\e \nm{x}\Ball\left(C(K^{(n_0)}\right).
\leqno{(3.4)}
$$
Setting $X_0=(R_{n_0}Q)^{-1} W$, we have from $(3.2)$, $(3.3)$, and
$(3.4)$ that
$$
U\left[(\nm{x}-\e)\Ball(Y_0)\right] \subset
X_0 + 2\e \nm{x}\Ball(X),
$$
which contradicts the Gorelik Principle since $\nm{x}$ is arbitrarily
large and $\e$ is arbitrarily small.\hfill\qed

\bs

\noindent{\bf 4. Uniform homeomorphisms between balls in ``close" spaces}

\ms

S. Mazur [Maz] proved that the unit balls of any two $\ell_p$, $1\le
p<\infty$, spaces are uniformly homeomorphic. This fact was extended by
E. Odell and T. Schlumprecht [OS], as a tool in solving the distortion
problem, to any two spaces with unconditional bases and nontrivial
cotype. It follows that any two balls in such spaces are uniformly
homeomorphic. However, the modulus of continuity of the uniform
homeomorphism or its inverse generally gets worse as the radius of
the balls  increases. In Corollary 4.7 we show that if the two spaces
are close to a common $\ell_p$ space then  ``large" balls are uniformly
homeomorphic with good modulus  (where ``large" depends on $p$ and how
``close" the two spaces are to  $\ell_p$.)

We begin by recalling the definition of the generalized Mazur map of
Odell and Schlumprecht. Let $ b=\{b_i\}_{i=1}^n\in\Rn{+n}$ with $\Vert
b\Vert_1=\sum b_i=1$ and let $X$ be a Banach space with a
$1$-unconditional basis
$\{e_i\}_{i=1}^\infty$. For
$x=\sum_{i=1}^\infty x_i e_i\in X$ with $x_i\ge0$, let
$G_X(b,x)=\prod_{i=1}^n x_i^{b_i}$, \ $G_X(b)=\sup_{\Vert
x\Vert=1}G_X(b,x)$
\ and $F_X(b)=F_{\ell_1,X}(b)$ the unique $x$ with the same support
as $b$ for which the
$\sup$ in the definition of $G_X(b)$ is attained.

Recall that, for  $x=\sum_{i=1}^\infty x_i e_i$ and $0<p<\infty$, we
denote $|x|^p=\sum_{i=1}^\infty |x_i|^p e_i$ and that the unconditional
basis
$\{e_i\}_{i=1}^\infty$ is said to be  $p$-convex  (respectively,
$q$-concave), with constant $C$, if $\Vert
(\sum_{n=1}^N |x^n|^p)^{1/p}\Vert\le C
(\sum_{n=1}^N \Vert x^n\Vert^p)^{1/p}$ (respectively, $\Vert
(\sum_{n=1}^N |x^n|^q)^{1/q}\Vert\ge C^{-1}
(\sum_{n=1}^N \Vert x^n\Vert^q)^{1/q}\ $) for all finite sequences
$\{x^n\}_{n=1}^N$  of vectors in $X$. (See [LT2,~section~1.d] for a
discussion of
$p$-convexity and
$q$-concavity).

Odell and
Schlumprecht proved that, if $X$ is $q$-concave with constant one
for some $q<\infty$ , then
$F_X$ extends naturally (i.e., homogeneously to the positive quadrant
in $\ell_1$, then extending to the other quadrants so that $F_X$ will
commute with changes of signs, and finally by continuity from $\ell_1^n$
to $\ell_1$)
 to a uniform homeomorphism of the unit balls of
$\ell_1$ and $X$.  Moreover, the modulus of continuity of $F_X$ and
$F_Y$ depend only on $q$.     If $X$ and
$Y$ both have $1$-unconditional basis and nontrivial cotype we shall
denote
$F_{X,Y}=F_{\ell_1,Y}\circ F^{-1}_{\ell_1,X}$. The following lemma
follows from a proof of Lozanovskii's theorem [Loz], (see, for
example, [T-J, Lemma~39.3], but note that misleading notation should
be corrected).

\proclaim Lemma 4.1.  Let $X$ be a Banach space with a
$1$-unconditional basis
$\{e_i\}_{i=1}^\infty$, let $b$ be a finitely supported vector in the
unit sphere of
$\ell_1$ with nonnegative coordinates,
 and set
$x=F_X(b)$. Then $f=\sum_{i=1}^\infty {{b_i}\over {x_i}}e_i^*$ is a
norming functional of $x$. Here $\{e_i^*\}_{i=1}^\infty$ denotes the
biorthogonal basis to $\{e_i\}_{i=1}^\infty$ and ${0\over 0}$ is
interpreted as $0$.

If $Y$ and $Z$ are two spaces with $1$-unconditional bases (or, more
generally, Banach lattices of functions) and $0<\theta<1$ we denote by
$X=Y^{\theta}Z^{1-\theta}$ the space of all sequences $x$ for which
$|x|$ can be written as $|x|=y^{\theta}z^{1-\theta}$ with $y\in Y,\ z\in
Z$ (the operations are defined pointwise) with the norm
$$
\Vert x \Vert=\inf\{\Vert y\Vert^{\theta}\Vert z\Vert^{1-\theta} : \
|x|=y^{\theta}z^{1-\theta}\}.
$$
By a theorem of Calderon [Cal], $Y^{\theta}Z^{1-\theta}$ is the
interpolation space, in the complex method, between $X$ and $Y$ with
parameter $1-\theta$ commonly denoted by $(Y,Z)^{1-\theta}$ (see, for
example, [BL]). The following lemma relates the functions
$F_X,\ F_Y$, and $F_Z$.

\proclaim Lemma 4.2. Let $b\ge 0$ have finite support with $\Vert
b\Vert_1=1$. Then for all spaces $Y$, $Z$ with $1$-unconditional bases,
$$
F_{Y^{\theta}Z^{1-\theta}}(b)=F_Y^\theta (b) F_Z^{1-\theta}(b).
$$

\Proof
$$\eqalign{
G_{Y^{\theta}Z^{1-\theta}}(b)=&\sup_{\Vert y\Vert^{\theta}\Vert
z\Vert^{1-\theta}=1}\prod y_i^{\theta b_i}\prod z_i^{(1-\theta)b_i}\cr
=&\sup_{ y\in Y}\left(\prod y_i^{\theta b_i}\sup_{ \Vert
z\Vert\sstwo{\scriptscriptstyle Z}=\Vert
y\Vert_{{\phantom {}}_{\scriptscriptstyle
Y}}^{-\theta/(1-\theta)}}\left(\prod
z_i^{b_i}\right)^{1-\theta}\right)\cr =&\sup_{y\in Y}\Vert
y\Vert_Y^{-\theta}\prod y_i^{\theta b_i} G_Z^{1-\theta}(b)\cr
=&G_Y^\theta (b) G_Z^{1-\theta}(b).
\cr}
$$
Let $y=F_Y(b),\ z=F_Z(b)$, and $x=y^\theta z^{1-\theta}$. Then
$\Vert x\Vert\le 1$ and
$$
G_{Y^{\theta}Z^{1-\theta}}(b,x)=\prod y_i^{\theta b_i}\prod
z_i^{(1-\theta)b_i}=G_Y^\theta (b)
G_Z^{1-\theta}(b)=G_{Y^{\theta}Z^{1-\theta}}(b).
$$
Consequently, $\Vert x\Vert=1$ and
$F_{Y^{\theta}Z^{1-\theta}}(b)=x$.\hfill
\qed

\ms

Recall that the $p$-convexification of a space $X$ with a
$1$-unconditional basis
(or a general Banach lattice of functions) is the space
$$
X^{(p)}=\{x=\sum x_ie_i \ :\ \Vert x \Vert_{\ss{X^{(p)}}}=\Vert
|x|^p\Vert^{1/p}<\infty\}.
$$
The expression $\Vert \cdot \Vert_{\ss{X^{(p)}}}$ is a norm whenever
$p\ge 1$. If
$p<1$, $X^{(p)}$ is sometimes  called the $1/p$-concavification of $X$.
It is still a Banach space if the original $X$ is $1/p$-convex with
constant one. (See [LT2, section 1.d] for more information on these
operations). Note that
$X^{(1/\theta)}=X^\theta\ell_\infty^{1-\theta}$.

\proclaim Corollary 4.3. If  $b$
is as in Lemma 4.2, then
$F_{X^{(1/\theta)}}(b)=F_X^\theta(b)$
for all spaces $X$ with  $1$-unconditional basis.

\Proof  Since $b$ has finite support, it is enough to prove the
corollary when $X$ is finite dimensional.  This case follows from
Lemma 4.2 and the the observation that for $b$ a strictly positive
vector in
$\ell_1^n$, \ $F_{\ell_\infty^n}(b)=(1,1,\dots,1)$.
\hfill\qed

\ms
For $1\le p \le \infty$, $p'$ denotes the conjugate index to $p$.

\proclaim Lemma 4.4. If $X$ is $p$-convex and $r\le
2p$ concave, both with constant one,  then
$X=\ell_p^{\theta}Y^{1-\theta}$ with $Y$ $p$-convex and
$2p$-concave, both with constant one, and
$\theta=1-{{2\over {(r/p){'}}}}={{2p}\over r}-1$.

\Proof Even if we just assume that $X$ is a space of functions on
the natural numbers, the convexity and concavity assumptions guarantee
that the unit vector basis is a  $1$-unconditional basis for $X$.
This implies that the general case follows from the case where $X$ is
finite dimensional, which we assume in the sequel to avoid irrelevant
topological and duality problems that arise in the infinite dimensional
setting.

Assume first that $p=1$. Since for any space $W$ with
$1$-unconditional basis,
$W^{(p)}=W^{1/p}\ell_\infty^{1-1/p}$,
$X^*=(X^{*({2\over{r{'}}})})^
{{2\over{r{'}}}}\ell_\infty^{1-{2\over{r{'}}}}$ so
$X=(X^{*({2\over{r{'}}})*})^
{1-\theta}\ell_1^{\theta}$ and $X^{*({2\over{r{'}}})*}$ is
$2$-concave with constant one.

For a general $p$, $X=(X^{({1\over p})})^
{1\over p}\ell_\infty^{1-{1\over p}}=(Z^
{1-\tau}\ell_1^{\tau})^{1\over p}\ell_\infty^{1-{1\over p}}$ \ \ with
$\tau=1-{2\over {(r/p){'}}}$ and $Z$ $2$-concave with constant one. This
can be rewritten as
$$
X=(\ell_1^{{1\over p}}\ell_\infty^{1-{{1\over p}}})^{\tau}
(Z^{{1\over p}}\ell_\infty^{1-{{1\over p}}})^{1-\tau}=
\ell_p^{\tau}(Z^{(p)})^{1-\tau}.
$$
\hfill\qed

We are now ready to investigate the modulus of continuity of
$F_{\ell_p,X}$ for spaces $X$ close to $\ell_p$.

\proclaim Proposition 4.5. Let $X$ be $p$-convex and $r$-concave, both
with constant one, and with
$r\le 2p$. Then the modulus of uniform continuity $\varphi(\epsilon)$ of
$F_{\ell_p,X}$ on the unit ball of $\ell_p$ is bounded by
$2(\varepsilon^{{2p\over r} -1}+(1-{p\over r})
\varphi_0(\varepsilon))$ for some function
$\varphi_0$ depending only on $p$.

\Proof Let $u,v$ be two nonnegative vectors of norm one in $\ell_p$. It
follows from Lemma 4.2 that, with the notation of Lemma 4.4,
$$
\leqalignno{
\Vert F_{\ell_p,X}(u)-&F_{\ell_p,X}(v)\Vert\sstwo{X}=\Vert u^\theta
F^{1-\theta}_{\ell_p,Y}(u)-v^\theta
F^{1-\theta}_{\ell_p,Y}(v)\Vert\sstwo{X}\cr
\le&\left\Vert w^\theta\left|
F^{1-\theta}_{\ell_p,Y}(u)-F^{1-\theta}_{\ell_p,Y}(v)
\right|\right\Vert\sstwo{X}+
\left\Vert |u^\theta-v^\theta|
F^{1-\theta}_{\ell_p,Y}(u)\vee
F^{1-\theta}_{\ell_p,Y}(v)\right\Vert\sstwo{X},&{(4.1)}\cr}
$$
where $w=\{w_i\}$ is defined by $w_i=u_i$ whenever the i-th component of
$F_{\ell_p,Y}(u)$ is smaller than that of $F_{\ell_p,Y}(v)$ and $w_i=v_i$
otherwise.

Now, the second term in 4.1 is dominated by
$$
\Vert |u-v|^\theta
F^{1-\theta}_{\ell_p,Y}(u)\vee F^{1-\theta}_{\ell_p,Y}(v)\Vert\sstwo{X}\le
\Vert u-v\Vert_p^\theta \Vert F_{\ell_p,Y}(u)\vee
F_{\ell_p,Y}(v)\Vert^{1-\theta}\sstwo{Y} \le
2^{1-\theta}\Vert u-v\Vert_p^\theta,\leqno{(4.2)}
$$
while the first term is dominated by
$$
\eqalign
{&(1-\theta)\left\Vert \left({w\over {F_{\ell_p,Y}(u)\wedge
F_{\ell_p,Y}(v)}}\right)^\theta
\left|F_{\ell_p,Y}(u)-F_{\ell_p,Y}(v)\right|\right\Vert\sstwo{X}\cr
\le&(1-\theta)\left\Vert \left({w\over {F_{\ell_p,Y}(u)\wedge
F_{\ell_p,Y}(v)}}\right)
\left|F_{\ell_p,Y}(u)-F_{\ell_p,Y}(v)\right|\right\Vert_{p}^\theta
\left\Vert F_{\ell_p,Y}(u)-F_{\ell_p,Y}(v)\right\Vert\sstwo{Y}^{1-\theta}\cr
=&(1-\theta)\left\Vert \left({w\over {F_{\ell_p,Y}(u)\wedge
F_{\ell_p,Y}(v)}}\right)^p
\left|F_{\ell_p,Y}(u)-F_{\ell_p,Y}(v)\right|^p\right\Vert_{1}^{\theta/p}
\left\Vert F_{\ell_p,Y}(u)-F_{\ell_p,Y}(v)\right\Vert\sstwo{Y}^{1-\theta}.
\cr}
$$
By  Corollary 4.3, $F_{\ell_p,Y}(u)=F_{\ell_1,Y}(u^p)=
F^{1/p}_{\ell_1,Y^{(1/p)}}(u^p)$. Thus, by Lemma 4.1,
$\left({u\over {F_{\ell_p,Y}(u)}}\right)^p$ and
$\left({v\over {F_{\ell_p,Y}(v)}}\right)^p$ are norm one functionals on
$Y^{(1/p)}$ and the first term in (4.1) is dominated by
$$
\eqalign
{(1-\theta)&\left\Vert \left({u\over
{F_{\ell_p,Y}(u)}}\vee{v\over {F_{\ell_p,Y}(v)}}\right)^p
\left|F_{\ell_p,Y}(u)-F_{\ell_p,Y}(v)\right|^p\right\Vert_{1}
^{\theta/p}
\left\Vert
F_{\ell_p,Y}(u)-F_{\ell_p,Y}(v)\right\Vert\sstwo{Y}^{1-\theta}\cr
\le&(1-\theta)2^{\theta/p}\left\Vert
\left|F_{\ell_p,Y}(u)-F_{\ell_p,Y}(v)\right|^p\right\Vert\sstwo{Y^{(1/p)}}
^{\theta/p}
\left\Vert
F_{\ell_p,Y}(u)-F_{\ell_p,Y}(v)\right\Vert\sstwo{Y}^{1-\theta}\cr
=&(1-\theta)2^{\theta/p}\left\Vert
F_{\ell_p,Y}(u)-F_{\ell_p,Y}(v)\right\Vert\sstwo{Y}.
\cr}
$$
Combining this with (4.1) and (4.2) and using the Odell and Schlumprecht
result we get that
$$\eqalign {
\Vert F_{\ell_p,X}(u)-F_{\ell_p,X}(v)\Vert\sstwo{X}\le&
2^{1-\theta}\Vert u-v\Vert_{p}^\theta+
(1-\theta)\varphi_0(\Vert u-v\Vert_{p})\cr
\le&2(\Vert u-v\Vert_{p}^{{2p\over r} -1}+(1-{p\over r})
\varphi_0(\Vert u-v\Vert_{p})\cr}
$$
for some function $\varphi_0$ depending only on $p$. \hfill\qed

\ms

Next we investigate the modulus of uniform continuity of $F_{X,\ell_p}$
for spaces close to $\ell_p$.

\proclaim Proposition 4.6. Let $X$ be $p$-convex and $r$-concave, both with
constant one, and with
$r\le 2p$. Then the modulus of uniform continuity $\varphi(\epsilon)$ of
$F_{\ell_p,X}^{-1}=
F_{X,\ell_p}$ on the unit ball of $X$ is bounded by \
$2(\epsilon+(1-{p\over r})^{1/p}\varphi_0(\epsilon))$ \ for some function
$\varphi_0$ depending only on $p$.

\Proof Write
$X^{(1/p)}=\ell_1^{\theta}Y^{1-\theta}$ with $Y$ $2$-concave
and
$\theta={{2p}\over r}-1$. Note that for $x\in X$ with $\Vert x\Vert=1$ and
$x\ge 0$,
$F_{\ell_p,X}^{-1}(x)= F_{X,\ell_p}(x)=(x\cdot f_x)^{1/p}$ where
$f_x$ is the norming functional of $x$. Let $x,x{'}$ be two norm one
vectors in the positive cone of $X$.
Writing $x^p=y\cdot f_y^{\theta}$,
${{f_x}\over{x^{p-1}}}=f_y^{1-\theta}$, we get from Lemma 4.2 that $y$
has norm one in $Y$ and $f_y$ norms it. Similarly  let $y{'}$ be the
norm one vector in $Y$ corresponding to $x{'}$. Let $z,z{'}$ be defined
by
$$
z(i)=\cases {x(i) \ \ if &$f_y(i)\le f_{y{'}}(i) $\cr x{'}(i)
&otherwise}
$$
and
$$
z{'}(i)=\cases {x{'}(i) \ \ if &$z(i)=x(i)$ \cr
x(i) &otherwise}\ \ \ .
$$
Then
$$\eqalign{
\Vert (x\cdot f_x)^{1/p}-&(x{'}\cdot f_{x{'}})^{1/p}\Vert_p=
\Vert |x\cdot ({{f_x}\over{x^{p-1}}})^{1/p}-x{'}\cdot
({{f_{x{'}}}\over{x^{{'}(p-1)}}})^{1/p}|^p\Vert_1^{1/p}\cr
\le&
2^{(p-1)/p}[\Vert
|x-x{'}|^p\cdot {{f_z{'}}\over{z^{{'}(p-1)}}}\Vert_1^{1/p} +
\Vert
z^p\cdot|({{f_x}\over{x^{p-1}}})^{1/p}-
({{f_{x{'}}}\over{x^{{'}(p-1)}}})^{1/p}|^p\Vert_1^{1/p}]\cr
\equiv &2^{(p-1)/p}[I+II].\cr}
$$
Now,
$$
I\le \Vert |x-x{'}|^p\Vert_{X^{(1/p)}}^{1/p}\Vert
\max\{{{f_x}\over{x^{p-1}}},{{f_{x{'}}}\over{x^{{'}(p-1)}}}\}
\Vert_{X^{{(1/p)}*}}^{1/p}\le 2^{1/p}\Vert x-x{'}\Vert_X
$$
by Lemma 4.1, and
$$\eqalign{
II\le &\Vert z^p\cdot|{{f_x}\over{x^{p-1}}}-
{{f_{x{'}}}\over{x^{{'}(p-1)}}}|\Vert_1^{1/p}\cr
\le &
\Vert z^p\cdot |f_y^{1-\theta}-f_{y{'}}^{1-\theta}|\Vert_1^{1/p}\cr
\le &(1-\theta)^{1/p}\Vert z^p\cdot f_y^{-\theta}\cdot
|f_y-f_{y{'}}|\Vert_1^{1/p}\cr
\le &(1-\theta)^{1/p}\Vert \max\{y,y{'}\}\cdot
|f_y-f_{y{'}}|\Vert_1^{1/p}\cr
\le &(1-\theta)^{1/p}[\Vert yf_y-y{'}f_{y{'}}\Vert_1+
\Vert |y-y{'}|\max\{f_y,f_{y{'}}\}\Vert_1]^{1/p}\cr
\le &(1-\theta)^{1/p}[\phi(\Vert y-y{'}\Vert)+2\Vert
y-y{'}\Vert]^{1/p}\cr
 }
$$
where $\phi$ is the uniformity function of $F_{Y,\ell_1}$ on the unit
ball of $Y$. Finally $\Vert
y-y{'}\Vert\le\psi(\Vert x-x{'}\Vert)$ where $\psi$ is the uniformity
function of $F_{X,Y}$. Note that $\phi$ and $\psi$ can be bounded by
functions depending only on $p$.\hfill\qed

\proclaim Corollary 4.7. Let $1\le p<\infty$. Then there is a function
$\varphi:\Rn{+}\rightarrow\Rn{+}$ such that for all
$1\le R <\infty$  there exists an $\varepsilon>0$ such that, if $X$ is
$p$-convex and
$(p+\varepsilon)$-concave, both with constant one, then there is a map
$\psi$ from the ball of radius
$R$ in $\ell_p$ onto the ball of radius $R$ in $X$ with the modulus of
continuity of $\psi$ and $\psi^{-1}$ bounded by $\varphi$.

\Proof Let $\psi$ be the natural (homogeneous) extension of
$F_{\ell_p,X}$ and use Propositions 4.5 and 4.6.  \hfill\qed

\bs

\noindent{\bf 5. Uniform homeomorphisms between $X$ and
$\ell_p\oplus X$}

\ms

In Corollary 5.3 we give a useful sufficient condition for
a space
$X$ to be uniformly homeomorphic to $X\oplus \ell_p$.
This is used in Corollary 5.7 to prove that $\Tp$ is
uniformly homeomorphic to $\Tp\oplus \ell_p$.

\proclaim Theorem 5.1. Let $1\le  p<\infty$. Let
$X=\sum_{n=0}^\infty X_n$ be a 1-unconditional sum of
spaces with 1-unconditional bases such that there are
norm preserving, homogeneous homeomorphisms
$f_n:\ell_p \ {\buildrel onto\over\longrightarrow}\  X_n$
with the property that, for some sequence
$R_n\rightarrow\infty$, the sequence
$f_{n_{|B(R_n)}}, f^{-1}_{n_{|B(R_n)}}$, $n=0,1,\dots$,
is equi-uniformly continuous. ( $B_Y(R)$, or just $B(R)$ if
there is no ambiguity, denotes the closed ball of radius
$R$ in $Y$.)   Then
$\ell_p\oplus X$ is uniformly homeomorphic to $X$.

For $X=(\sum\oplus\ell_{p_n})_q$ and $p_n\to p$, this theorem
was proved by M. Ribe [Rib3] for $p=1$ and extended by I.
Aharoni and J. Lindenstrauss [AL2] for $p>1$. The proof we
sketch closely follows a simplification of the proof in [AL2]
given by Y. Benyamini in his nice exposition [Ben2].

In what follows we shall frequently use the spaces
$\ell_p\oplus X_n$, $\ell_p\oplus X_n\oplus X_{n+1}$ and
other spaces of similar nature. Since the isometric nature
of the spaces we deal with is important in the proof, we
emphasize that by, e.g.,
$\ell_p\oplus X_n\oplus X_{n+1}$ we mean the space of all
triples $(u,x_n,x_{n+1})$, $u\in \ell_p,\ x_i\in
X_i,\ i=n,n+1$, with norm
$$
\Vert (u,x_n,x_{n+1})\Vert=(\Vert u\Vert^p+\Vert
(0,\dots,0, x_n,x_{n+1},0\dots)\Vert_X^p)^{1/p}.
$$
By coupling the ${X_n}'s$, we may assume that for each
$n$, $X_n =Y_n\oplus Z_n$ and there are norm
preserving homogeneous homeomorphisms
$g_n:\ell_p\ {\buildrel onto\over\longrightarrow}\  Y_n$
and $h_n:X_n\ {\buildrel onto\over\longrightarrow}\  Z_n$
such that the sequence
$g_{n_{|B(R_n)}}, g^{-1}_{n_{|B(R_n)}}, h_{n_{|B(R_n)}},
h^{-1}_{n_{|B(R_n)}}$,
$n=0,1,\dots$, is equi-uniformly continuous.

Define now, $I_n:\ell_p\oplus X_n\ {\buildrel
onto\over\longrightarrow}\  X_n$ by
$$
I_n(u,x)=I_n(u,(y,z))={{\Vert
(u,x)\Vert}\over{\Vert(g_n(u),h_n(x))\Vert}}
(g_n(u),h_n(x)).
$$
Clearly, the ${I_n}'s$ are norm preserving, homogeneous,
and
$I_{n_{|B(S_n)}},I^{-1}_{n_{|B(S_n)}}$, $n=0,1,\dots,$ form
an equi-uniformly continuous family for some
sequence $S_n\to\infty$. Since, in order to prove the
theorem we may pass to a subsequence of the ${X_n}'s$,
we may assume that $S_n>2^{n+1}$.

The following proposition imitates Proposition 6.2 in
[Ben2].

\proclaim Proposition 5.2. There is a family of
homeomorphisms, $\{F_t\ ;\ 0\le t<\infty\}$, so that
\item {(i)} For $n=0,1,\dots$ and $2^{n-1}\le t< 2^n,\
(0\le t< 1\ for\  n=0) $, $F_t$ maps $\ell_p\oplus
X_n\oplus X_{n+1}$ onto  $X_n\oplus X_{n+1}$.
\item{(ii)}
$F_t$ is homogeneous and norm preserving, $0\le
t<\infty$.
\item {(iii)}
$F_{2^{n-1}}(x,y,z)=(I_n(x,y),z)$,
$n=1,2,\dots$
\item {(iv)} $F_{2^{n}}(x,y,z)=(y,I_{n+1}(x,z))$,
$n=0,1,\dots$
\item{(v)} Denote, for $2^{n-1}\le t<2^n$,
$G_t=F_{t_{|B_{\ell_p\oplus
X_n\oplus X_{n+1}}(t)}}$. Then the families
$\{G_t(a)\}$, $\{G_t^{-1}(a)\}$ are equi-uniformly continuous
in both $t$ and $a$.

Before we sketch the proof of the proposition let us
apply it to give the

\noindent{\bf Proof of Theorem 5.1.} Define an
homeomorphism $f:\ell_p\oplus X\ {\buildrel
onto\over\longrightarrow}\  X$ in the following way: For
$a=(u,x_0,x_1,\dots)\in \ell_p\oplus X$ with $2^{n-1}\le
\Vert a\Vert<2^n$ ($0\le \Vert a\Vert<1$ for $n=0$) let
$$
f(a)=(x_0,x_1,\dots,x_{n-1},F_{\Vert a
\Vert}(u,x_n,x_{n+1}),x_{n+2},\dots)
$$
and
$$
g(a)=f(a)/\Vert f(a)\Vert.
$$
For $(x_n,x_{n+1})\in X_n\oplus X_{n+1}$ and $2^{n-1}\le
t<2^n$ denote
$(u^t,y^t_n,y^t_{n+1})=F^{-1}_t(x_n,x_{n+1})$. Then it is
easily checked that for $x=(x_0,x_1,\dots)$
$$
g^{-1}(x)={{(u^{\Vert x\Vert},x_0,x_1,\dots,x_{n-1},
y^{\Vert x\Vert}_n,y^{\Vert
x\Vert}_{n+1},x_{n+2},\dots)}\over
{\Vert (u^{\Vert x\Vert},x_0,x_1,\dots,x_{n-1},
y^{\Vert x\Vert}_n,y^{\Vert
x\Vert}_{n+1},x_{n+2},\dots)\Vert}}.
$$
The 1-unconditionality of the decomposition of $X$
to $\sum\oplus X_n$ and Proposition 5.2 now show that $g$
and $g^{-1}$ are uniformly continuous.\hfill\qed

\ms

\noindent{\bf Proof of Proposition 5.2.} Consider the
maps
$$
\varphi_n:\ell_p\oplus X_n\oplus X_{n+1}\to
\ell_p\oplus_p\ell_p\oplus_p\ell_p
$$
and
$$
\psi_n:X_n\oplus X_{n+1}\to
\ell_p\oplus_p\ell_p
$$
given by
$$
\varphi_n(u,x,y)={{\Vert
(u,x,y)\Vert}\over{\Vert(u,f_n(x),f_{n+1}(y))\Vert}}
(u,f_n(x),f_{n+1}(y))
$$
and
$$
\psi_n(x,y)={{\Vert
(x,y)\Vert}\over{\Vert(f_n(x),f_{n+1}(y))\Vert}}
(f_n(x),f_{n+1}(y)).
$$
Note that these are norm preserving homogeneous maps and
the sequence \hfill\break
\hbox{$\varphi_{n_{|B(2^{n+1})}},\psi_{n_{|B(2^{n+1})}},
\varphi^{-1}_{n_{|B(2^{n+1})}},\psi^{-1}_{n_{|B(2^{n+1})}}
$}
is equi-uniformly continuous. Assume we can prove the
proposition in the special case where
$X_n=\ell_p$ for all
$n$ and
$\sum\oplus X_n=
\sum\oplus_p X_n$ and call the resulting maps $H_t$
instead of $F_t$. Then for
$n=0,1,\dots$ and $2^{n-1}\le t<2^n$
define $F_t=\psi^{-1}_n\circ H_t\circ
\varphi_n$. Clearly $F_t$ satisfy the conclusions of the
proposition.

The construction of $H_t$ is given in steps II-IV of the
proof of Proposition 6.2 of [Ben2]. The only
additional thing to notice is that all the maps there are
homogeneous.\hfill\qed

\proclaim Corollary 5.3. Let $1\le p<\infty$ and
let
$X=\sum X_n$ be an unconditional sum of spaces $X_n$ with
unconditional bases which are
$p$-convex with constant $C$ and
$(p+\varepsilon_n)$-concave with constant $C$, for some
positive sequence
$\varepsilon_n\to 0$ and some $1\le C<\infty$. Then
$\ell_p\oplus X$ is uniformly homeomorphic to $X$.

\Proof We may assume that the unconditionality constant
of the sum $\sum X_n$ is $1$. By a theorem from [FJ], each
$X_n$ is isomorphic, with constant independent of $n$, to a
space with a 1-unconditional basis which is $p$-convex  and
$(p+\varepsilon_n)$-concave with constants one. Applying
 Corollary 4.7 stated after Proposition 4.6 (and
the above isomorphisms) we get the maps $f_n:\ell_p\to
X_n$ as in the statement of  Theorem 5.1. Now apply
Theorem 5.1. \hfill\qed

\ms

We shall denote by
$\Tp$ (or just $\T$ when $p=1$) the closed span of a certain
subsequence of the unit vector basis for the
$p$-convexification of the modified Tsirelson space.  In
fact, it is known (see [CO], [CS]) that this space is (up to
an equivalent renorming) the space which is usually denoted by
$\Tp$; that is, the $p$-convexification of  the dual to
Tsirelson's original space; but we avoid using this for the
convenience of the reader.

 The construction of
$\T$, as well as an elementary proof that
$\T$ is reflexive, is contained in [Joh2].  The space
$\T$ is the completion of $c\ss{{00}}$ under the norm defined
implicitly  by the formula:

$$  ||x|| = \max \left\{ ||x||_{c_0}, \ {1\over 2}
\sup \{ \sum_{j=1}^{(n+1)^n} ||E_j
x||\}\right\},\leqno{(5.1)}
$$ where the $\sup$ is over all  sequences of disjoint finite
sets $E_j$ of positive integers  with
$j\le E_j$ for each $1\le j\le n$. Such sequences are
called {\sl allowable}.  The choice of the growth $n\mapsto
(n+1)^n$ is simply one of convenience; the results of [CJT],
[CO] imply that the space is the same if e.~g. ``$n$"
instead of ``$(n+1)^n$" is used.

In order to show that for all $1\le p < \infty$, $\Tp$,
the
$p$-convexification of $\T$, is uniformly homeomorphic to
$\Tp\oplus\ell_p$, we need that for any $q>1$ the span of
a sufficiently thin subsequence of the unit vector basis
$\en$ for $\T$ has $q$-concavity constant independent of $q$.
(In fact, the subsequence can be just a ``tail" of the
basis.) The proof uses an elementary lemma:

\proclaim Lemma 5.4. Suppose $q>1$, $\e>0$, and
$1=m_1<m_2<\dots$ with
$\sum_{k=2}^\infty (m_k-m_{k-1})^{1-q} <\e$.  Then for all
sequences
$a_1\ge a_2 \ge \dots \ge 0$,
$$
\sup_k \sum_{j=m_k}^{m_{k+1}-1} a_j \ge
(1+\e)^{-1/q}\left(\sum_j a_j^q\right)^{1/q}.
$$

\Proof  Assume, without loss of generality, that the left side
is one.  Then

$$
\eqalign{
\sum_j a_j^q & = \sum_{k=1}^\infty \sum_{j=m_k}^{m_{k+1}-1}
a_j^q
 \le \sum_{k=1}^\infty \left( \sum_{j=m_k}^{m_{k+1}-1}
a_j\right) a_{m_k}^{q-1} \le
\cr
\left(\sup_k \sum_{j=m_k}^{m_{k+1}-1} a_j\right)
\left(\sum_{k=1}^\infty a_{m_k}^{q-1}\right) &\le
\left(\sup_k \sum_{j=m_k}^{m_{k+1}-1} a_j\right)
\left(a_{m_1}^{q-1}+\sum_{k=2}^\infty
(m_k-m_{k-1})^{1-q}\right), \cr}
$$ so the desired conclusion follows.\hfill\qed

\proclaim Corollary 5.5. There is a constant $M<\infty$ so
that for any
$1<q<{3\over 2}$ the span of $\{e_n\}_{n=N(q)}^\infty$ in
$\T$ has
$q$-concavity constant less than $M$, where
$N(q)=2\superscript{{q\over{q-1}}}$.

{\bf Proof:} A general result of Maurey's [Mau2] (see [LT2,
proof of Theorem 1.f.7]) shows that the  $q$-concavity
constant of a  Banach lattice is at most a constant multiple
of the lower
$s$ constant of the lattice, where $q^*=s^*+1$; i.e.,
$s={1\over {2-q}}$. (Trace the constants in the proof in
[LT2] for the dual statement in the special case $2 < r=p+1$.)

Thus it is enough to prove the corollary with  ``$q$-concavity
constant" replaced with ``lower $q$ constant"  and ``$N(q)$"
replaced with ``$M(q)\equiv N(s)$" where, as above,
$s={1\over {2-q}}$, so that
$M(q)=2\superscript{{1\over{q-1}}}$.

So let $1<q<{3\over 2}$ be fixed and suppose that
$\seq{x}{n}$ are disjoint vectors in $\T\cap c\ss{{00}}$
which are supported on
$[M(q),\infty)$.  We can assume that $||x_1||\ge ||x_2||\ge
\dots ||x_n||$.  If, for some $m$, $F$ is a set of at most
$(m+1)^m$ indices so that for each $j$ in $F$,  $x_j$ is
supported in $[m,\infty)$, then by the definition of $\T$,
$\displaystyle ||\sum_{j=1}^n x_j|| \ge ||\sum_{j\in F} x_j||
\ge  {1\over 2} \sum_{j\in F} ||x_j||$.

Define $m_1=1$ and $m_{k}=m_{k-1}+M(q)^{k-1}$ for $k>1$.
Since for $k>1$, $m_{k+1} - m_k \le (m_k + 1)^{m_k}$ and
since at most $m_k$ of the $x_i$'s, $1\le i \le
m_{k+1}$, fail to be supported on $[m_k,\infty)$, we get
$$ ||\sum_{j=1}^n x_j|| \ge  {1\over 2} \sup_k
\sum_{j=m_k}^{m_{k+1}-1} ||x_j||.
$$ The sequence $\mk$ satisfies the condition of the lemma
with
$\e=1$, and thus the lemma yields that that
$\{e_n\}_{n=M(q)}^\infty$ has a lower $q$ estimate with
constant
$2^{-1-{1/q}}$.\hfill\qed

\ms

If $X$ has a $1$-unconditional basis and is $q$-concave with
constant $M$, then its $p$-convexification is $pq$-concave
with constant $M^{1\over p}$. Thus the next proposition is
an immediate consequence of Corollary 5.5.

\proclaim Proposition 5.6. There is a constant $M<\infty$
so that for any $1\le p <r < {3p\over 2}$ the span of
$\{e_n\}_{n=N(p,r)}^\infty$ in $\Tp$ has $r$-concavity
constant less than $M^{1/p}$, where
$N(p,r)=2\superscript{{r\over{r-p}}}$.

As an immediate consequence of Proposition 5.6 and Corollary
5.3 we get:

\proclaim Proposition 5.7. For $1\le p < \infty$, $\Tp$ is
uniformly homeomorphic to $\Tp\oplus\ell_p$.

Finally, by combining Proposition 5.7 and Proposition 2.8 we get:

\proclaim Theorem 5.8.  Let $1<p\ss{1}<p\ss{2}<\dots<p\ss{n}
<2$ or  $2<p\ss{1}<p\ss{2}<\dots<p\ss{n}  <\infty$.   Set
$\displaystyle X=\sum_{k=1}^n \Tpk$.  If
$X$ is uniformly homeomorphic to a Banach space $Y$ if and
only if   $Y$ is isomorphic to $X\oplus
\sum_{k\in F} \lpk$ for some subset $F$ of $\{1,2,\dots,n\}$.
Consequently, $X$ is uniformly homeomorphic to exactly $2^n$
mutually nonisomorphic Banach spaces.

\Proof $X$ is uniformly homeomorphic to spaces $Y$ of the
desired form by Proposition 5.7.  As mentioned in section 2,
the (obvious) fact that every $n$-tuple (even
$(n+1)^n$-tuple) of disjoint unit vectors in $\Tp$ whose
first $n$ coordinates are zero is $2$-equivalent to the
unit vector basis of $\lpn$ implies that $\Tp$ is as.
$\sLp$.  Of course, $\ell_r$ does not embed into $\Tp$ for any
$r$ (the general case is immediate from the special case
$p=1$, proved  in [Joh2]).  So it remains to observe that
$\Tp$ is isomorphic to its hyperplanes.  In fact, the
right shift is an isomorphism on $\T$, hence on $\Tp$ for
all $1\le p < \infty$.  This is included in the first result
in [CJT] for the classical Tsirelson space. To avoid
using the fact that this space is the same space (up to an
equivalent renorming) that we denote by $\T$, the reader will
have to check directly that the right shift is an isomorphism
on $\T$.   So Proposition 2.8 applies to show that any
uniform homeomorph of $X$ must be of the form for $Y$ in the
statement of Theorem 5.8.  Finally, $\lpj$ embeds into
$X\oplus
\sum_{k\in F} \lpk$ iff and only if $j$ is in $F$, so all
the $Y$'s of the mentioned form are isomorphically distinct.
\hfill\qed

\bs

\noindent{\bf 6. The nonseparable case}

\ms

Some of the uniqueness theorems proved in earlier sections
have nonseparable analogues.  In particular, there is a
nonseparable version of Theorem 2.1:

\proclaim Theorem 6.1. If $X$ is a Banach space which is
uniformly homeomorphic to $\ell_p(\Gamma)$, $1<p<\infty$,
then
$X$ is isomorphic to $\ell_p(\Gamma)$.

\Proof  A back-and-forth argument yields that each
separable subspace of $X$ is contained in another subspace
of $X$ which is uniformly homeomorphic to $\ell_p$, hence
isomorphic to $\ell_p$. It is then easy to see that there
is a $\lam<\infty$ so that each separable subspace of $X$
is contained in another subspace of $X$ which is
$\lam$-isomorphic to $\ell_p$. So $X$ is $\sLp$ and hence
is isomorphic to a complemented subspace of $L_p(\mu)$ for
some measure $\mu$ [LPe], and of course $\ell_2$ does not
embed into $X$. The conclusion then follows from the
next lemma, which is implicit in [JO] although not
stated there.
\ms

\proclaim Lemma 6.2. Let $X$ be  an infinite dimensional
complemented subspace of $L_p(\mu)$ for
some measure $\mu$,  $1<p<\infty$, and suppose that no
subspace of $X$ is isomorphic to $\ell_2$. Then
$X$ is isomorphic to
$\ell_p(\Gamma)$, where $\Gamma$ is the density character
of $X$.

\Proof Since every Hilbertian subspace of $L_p$ is
complemented for $2<p<\infty$ [KP], by duality we can
assume that $2<p<\infty$.  By Theorem 3 of [JO], $X$
embeds into $\ell_p(\Gamma)$, and, as noted in the proof
of that theorem, $X$ is isomorphic to a space of the form
$\displaystyle
\left(\sum_{\gamma\in\Gamma}X_\gamma
\right)_{\ell_p(\Gamma)}$
with each space $X_\gamma$ separable (and, without loss
of generality, infinite dimensional). But then the
$X_\gamma$'s are uniformly isomorphic to uniformly
complemented subspaces of $\ell_p(\Gamma)$, hence by [JO]
and [JZ] (or see [Joh1]) are uniformly isomorphic to
$\ell_p$.
\hfill\qed

\ms

The isomorphic classification of
$C(K)$ for compact metric spaces $K$ [BP] implies that
$c_0$ is isomorphic to $C(K)$ if and only if the
$\omega$-th derived set $K^{(\omega)}$ of $K$ is empty. In
[DGZ] it is proved that for any compact Hausdorff space
$K$, if
$K^{(\omega)}$ is empty, then $C(K)$ is uniformly
homeomorphic (even Lipschitz equivalent) to $c_0(\Gamma)$
with
$\Gamma$ the density character of $C(K)$.  Part of the
interest of this result is that, in contrast to the
separable case,  in the nonseparable setting this can
happen with
$C(K)$ not isomorphic to $c_0(\Gamma)$, [AL1].  Gilles Godefroy
pointed out to us that Corollary 3.2 and previously known
results yield a converse to the mentioned theorem from
[DGZ] and kindly suggested that we include the result in
this paper.

\proclaim Theorem 6.3. (G. Godefroy) Suppose that $K$ is
a compact Hausdorff space and $C(K)$ is uniformly
homeomorphic to $c_0(\Gamma)$.  Then
$K^{(\omega)}$ is empty.

\Proof The usual back-and-forth argument  shows that
every separable subspace of $C(K)$ is contained in a
subalgebra which is uniformly homeomorphic to $c_0$ hence
isomorphic to $c_0$  by Corollary 3.2. This
implies that $K$ is scattered
by [PS] (or see [Sem, Theorem 8.5.4]); that is,
whenever $K^{(\a)}$ is nonempty, $K^{(\a+1)}$ is a
proper subset of $K^{(\a)}$.  But if $K$ is scattered and
$K^{(\omega)}$ is nonempty, then $\omega^\omega$ is a
continuous image of $K$ by a result of Baker's [Bak].
But then $C(\omega^\omega)$ is isomorphic to a subspace of
$C(K)$ and hence,  by the first sentence of this proof,
$C(\omega^\omega)$ is isomorphic to a subspace of $c_0$,
which of course is false.
\hfill\qed

\ms
\noindent {\bf Remark.}  Theorem 6.3 and the mentioned
result in [DGZ] yield that if $c_0(\Gamma)$ is uniformly
homeomorphic to a $C(K)$ space $X$, then $c_0(\Gamma)$ is
Lipschitz equivalent to $X$.  We do not know whether this
holds for a general space $X$ even in the separable case.

\bs

\noindent{\bf 7. Spaces determined by their finite dimensional
subspaces}
\bs

 We begin
with some  notions and definitions.
\ms

 A {\sl paving\/} of a (separable) Banach space $X$ is a
sequence $E_1\subset E_2\subset\cdots$ of finite dimensional
subspaces of $X$ whose union is dense in $X$.
We say that the two  Banach spaces $X$ and $Y$  {\sl have a common
paving\/} if there exist pavings
$\En$ of $X$ and $\Fn$ of
$Y$  so that the isomorphism constants between $E_n$ and $F_n$
are uniformly bounded.

\ms

\noindent{\bf Example} (warning).  $c_0$ and $c_0\oplus\ell_1$ have a
common paving, since both are paved by
$\ell_{\infty}^{2^n}\oplus_\infty \ell_1^n$.

\ms

\noindent{\bf Definition 7.1.} A separable Banach space $X$ is
said to be {\sl determined by its pavings\/} provided that $X$ is
isomorphic to every separable Banach space which has a common
paving with
$X$.

\ms

\noindent{\bf Definition 7.2.}  $X$ and $Y$ are said to have the
{\sl same finite dimensional subspaces\/} if each is finitely
crudely representable in the other.
A separable Banach space $X$ is
said to be {\sl determined by its finite dimensional subspaces\/}
provided that
$X$ is isomorphic to every separable Banach space  which has the
same finite dimensional subspaces as $X$.

\ms

Obviously, every Banach space which is determined by its finite
dimensional subspaces is determined by its pavings.  It is easy
to see and well known that $\ell_2$ is determined by its
pavings.  In this section we shall show, among other things,
that any space determined by its pavings must be ``close" to
$\ell_2$, but that there are spaces not isomorphic to $\ell_2$
which are determined by their pavings.  On the other hand, we
make the following conjecture:

\ms

\noindent{\bf Conjecture  7.3.}   Every (separable,
infinite dimensional) Banach space which  is determined by its
finite dimensional subspaces is   isomorphic to $\ell_2$.

\ms

We start the discussion of the concepts introduced above by
proving:

\ms

\proclaim Proposition 7.4.  If $X$ is determined by its pavings,
then for every $\e>0$, $X$ has type $2-\e$ and cotype
$2+\e$.

\Proof   Otherwise $\ell_p$ is finitely representable in
$X$ for some $1\le p\neq 2 \le \infty$ by the
Krivine-Maurey-Pisier theorem ([Kri], [MP]).
 If $\Fn$ is a  paving of any subspace of $L_p$, then using
the fact that $\ell_p$ is finitely representable in every
finite codimensional subspace of $X$, it is easy to see that
$X$ has a paving of the form
$\{E_{n}\oplus_p F_{n}\}_{n=1}^\infty$ with $\En$ itself a
paving of $X$.  This means that $X\oplus Y$ has a common paving
with
$X$ for every subspace $Y$ of $L_p$, which implies that every
subspace of
$L_p$ is isomorphic to a complemented subspace of $X$. But in
[JS] it was proved that for $2<p<\infty$, no separable Banach
space is complementably universal for subspaces of $\ell_p$.
Probably  the same result is true when $1\le p<2$, although we can
get by with the weaker fact that there is no separable  Banach
space which is complementably universal for subspaces of $L_p$
for $p$ in this range.  This weaker fact follows from the proof
of what is called in [JS] the Basic Result and two facts which
can be found in [LT1,2]: (1)
$\ell_r$ embeds into $L_p$ when $1\le p < r <2$; and (2) for each
$r\neq 2$, $\ell_r$ has a subspace which fails the compact
approximation property.
\hfill\qed

Following [Joh3], we say that $X$ has property $H(m,n,K)$
provided that there exists an $m$-codimensional subspace of $X$,
all of whose
$n$-dimensional subspaces are $K$-isomorphic to $\ell_2^n$.
So a Banach space $X$ is {\sl asymptotically Hilbertian\/} in
Pisier's terminology [Pis] or as. $\sLt$ in our terminology if and
only if there exists a
$K$ so that for every $n$, $X$ has $H(m,n,K)$ for some $m=m(n)$.

\proclaim Lemma 7.5.   Given $m$, $k$, $K$, and a
separable Banach  space
$X$, the following are equivalent:
\item{(i)} $X$ has property $H(m,k,K)$.
\item{(ii)} There is a paving $\En$ for $X$ such that for every
$n$, $E_n$ has property $H(m,k,K)$.
\item{(iii)} Every finite dimensional subspace of $X$ has property
$H(m,k,K)$.
\item{(iv)} Every subspace of $X$ has property $H(m,k,K)$.
\item{(v)} Every space which is finitely representable in $X$
has property $H(m,k,K)$.

\Proof  The implications
(v)$\implies$(i)$\implies$(iv)$\implies$(iii)$\implies$(ii) are
all obvious.  To prove (ii)$\implies$(v), it is enough to
observe that ultraproducts of $\En$ have property $H(m,k,K)$
since any separable space finitely representable in $X$ embeds
isometrically into an ultraproduct of any paving of
$X$.\hfill\qed

\ms

Notice that the preceding lemma implies that
asymptotically Hilbertian is a local property.

The next proposition is suggested by a result and proof due to P.
Casazza [CS,~Theorem~Ae13].

\proclaim Proposition 7.6. Suppose that $X$ is asymptotically
Hilbertian, $\En$ is a paving for $X$, and $\tilde E$ is an
ultraproduct of  $\En$.  Then $\tilde E$ is isomorphic to
$X\oplus H$ for some (nonseparable) Hilbert space $H$.

\Proof The subspace of $\tilde E$ consisting of all Cauchy
sequences $\xn$ in $X$ with $x_n$ in $E_n$ is easily seen to be
isometrically isomorphic to $X$. $\tilde E$ is the direct
sum of this subspace and the subspace of all $\xn$ in $E_n$
for which $\xn$  weakly converges along $\U$ to zero; this
latter subspace is isomorphic to a Hilbert space.

\proclaim Corollary 7.7.  If $X$ is asymptotically Hilbertian
and
$Y$ has a common paving with $X$, then $X\oplus \ell_2$
is isomorphic to $Y\oplus \ell_2$.  Consequently, if also
every infinite dimensional subspace of $X$ contains a complemented
subspace which is isomorphic to $\ell_2$ (in particular, if
$X$ is
an $\ell_2$ sum of finite dimensional spaces), then $X$ is
determined by its pavings.

\Proof   If $\En$ is a paving of $X$ which is uniformly
equivalent to a paving of $Y$, and $\tilde E$ is an ultraproduct
of  $\En$, then we get from Proposition 7.6 that
$X\oplus H$ is isomorphic to $Y\oplus H$ for some nonseparable
Hilbert space $H$, from which the first conclusion follows
easily. For the ``consequently" statement, we see that $X$ is
isomorphic to $X\oplus \ell_2$ since $X$ contains a complemented
copy of $\ell_2$.  In fact, every subspace of $X\oplus \ell_2$
(in particular, $Y$)   contains a complemented copy of $\ell_2$,
so also $Y$ is isomorphic to $Y\oplus \ell_2$.

\ms

This corollary implies that there are spaces not isomorphic to
$\ell_2$ which are determined by their pavings.  In order to
continue our investigation, we introduce a further notion:

\ms

We say that $X$ is  {\sl finitely
homogeneous}
provided there is a $K$ so that $X$ is  finitely crudely
representable with constant $K$ in every finite codimensional
subspace of itself.

\ms

First we observe that finite homogeneity is determined by finite
dimensional subspaces.

\proclaim Proposition 7.8.  $X$ is finitely homogeneous if and
only if there is a $K$ so that for every finite dimensional
subspace
$E$ of $X$, there is a space $Z=Z(E)$ which has a monotonely
subsymmetric and unconditional Schauder decomposition into spaces
isometric to $E$ and which is  finitely crudely representable in
$X$ with constant $K$.

\Proof  The forward direction involves just a combination
of modern refinements of the Mazur technique  for
constructing basic sequences and the Brunel-Sucheston
construction; see e.g. [Pis, Chpt. 14]. The reverse direction is
even easier (if
$F$ is a finite subset of
$X^*$ and
$E\oplus E\oplus\cdots\oplus E$ is a ``good" Schauder
decomposition of some subspace of $X$ into $N$-copies of $E$ with
$N$ large relative to the cardinality of $F$, then the action of
$F$ on two different copies of $E$ is essentially the same.  Now
pass to  differences.)\hfill\qed

\proclaim Proposition 7.9. If the separable Banach space $X$ is
finitely homogeneous and
$X$ has the same finite dimensional subspaces as the separable
Banach  space $Y$, then $X$ and $Y$ have a common paving.

\Proof The space $Y$ is also finitely homogeneous by
Proposition 7.8.  Now it is a simple exercise to show that if
$\En$ is a paving for $X$ and $\Fn$ is a paving for $Y$, then
there is a  subsequence of  $\EnFn$ which is uniformly equivalent
to  pavings for both $X$ and $Y$. \hfill\qed

\ms

The next proposition crystalizes some of what has already been
used.

\proclaim Proposition 7.10. If there is a $K$ so that the
separable Banach space $Y$ is finitely crudely representable
with constant $K$ in every finite codimensional subspace of the
separable Banach space $X$, then $X\oplus Y$ and $X$ have a
common paving.

\proclaim Corollary 7.11. If $X$ is finitely homogeneous and
determined by its pavings, then every separable Banach
space which is finitely representable in $X$ is isomorphic
to a complemented subspace of $X$. Moreover, $X$ has an
unconditional Schauder decomposition such that the summands form
a  sequence which is dense in the set of all finite dimensional
subspaces of
$X$, and hence every separable Banach space which is finitely
representable in
$X$ has the bounded approximation property.  Finally, if $X$ has
GL-l.u.st.; respectively, GL;  then so does every subspace.

See [T-J, \S34] for the definitions of GL-l.u.st. and GL
properties.

Using James' characterization of reflexivity/superreflexivity
[Jam] and the limiting technique of Brunel-Sucheston
(again as in [Pis, Chpt. 14]), we get from Proposition 7.13:

\proclaim Proposition 7.12.  If $X$ is reflexive and is
determined by its pavings, then $X$ is superreflexive.

We do not know whether there exist nonreflexive spaces which are
determined by their pavings.

Proposition 7.8 implies that if $X$ is finitely homogeneous,
then
$X\oplus X$ is finitely crudely representable in $X$.  The
properties are not equivalent, however, since every Tsirelson type
space is isomorphic to its square. The weaker condition gives
similar conclusions to those in Corollary 7.12 for spaces which
are determined by their finite dimensional subspaces:

\proclaim Proposition 7.13. If $X$ is determined by its finite
dimensional subspaces and $X\oplus X$  is finitely crudely
representable in $X$, then $X$ is isomorphic to
$X\oplus Y$ for every space $Y$ which is finitely representable
in $X$.  Hence if $X$ has GL-lust, then so does every subspace.

Corollary 7.7, Proposition 7.13 and [KT] yield:

\proclaim Proposition 7.14. Let $k_n\rightarrow\infty$
appropriately, set
$E_n=\span \{e_j\}_{j=k_n}^{k_{n+1}-1}$ in $\Tt$, and let
$X$ be the $\ell_2$-sum of $\En$.  Then $X$ is determined by
its pavings but $X$ is not determined by its finite
dimensional subspaces.

\Proof  Corollary 7.7 says that $X$ is determined by its
pavings no matter how $\kn$ is defined.  Next note that
$E_n\oplus E_n$ is uniformly isomorphic to $E_n\oplus
\ell_2^{k_{n+1}-k_n}$; again, no matter how $\kn$ is defined.
This follows from the equivalence of
$\eevenn$ and $\eoddn$ with $\en$ and the fact that the last half
of a finite string of $e_j$'s is equivalent to an orthonormal
sequence.  This yields that $X$ is isomorphic to its square.
Thus by Proposition 7.12, if $X$ were determined by its finite
dimensional subspaces, then every subspace of $X$ would have
GL-l.u.st., which easily implies that there is a $K$ so that
every finite dimensional subspace of $X$ has GL-l.u.st. constant
less than $K$.  But by [KT],
 $\Tt$ has a subspace which fails GL-l.u.st., so for each
$N$, $\span\{e_j\}_{j=N}^\infty$ has finite dimensional
subspaces with arbitrarily large GL-l.u.st. constant.
Consequently, having defined $k_n$, $X_n$ will have a subspace
with GL-l.u.st. constant larger than $n$ if $k_{n+1}$ is
sufficiently large.
\hfill\qed

\proclaim Proposition 7.15. Suppose that $p_n\downarrow 2$,
$k_n\rightarrow \infty$, and let $X$ be the $\ell_2$-sum of
$\ell_{p_n}^{k_n}$.  If $X$ is determined by its finite
dimensional subspaces, then $X$ is isomorphic to $\ell_2$.

\Proof    Assume without loss of generality that all the
$p_n$'s are smaller than $4$.  Now if
$2<p<4$, then $\ell_p^k$ is uniformly isomorphic to
$\ell_p^{k-\sqrt k}\oplus \ell_2^{\sqrt k}$ ([BDGJN], [JSc],
which yields that
$X$ is isomorphic to $X\oplus W$, where $W$ is the $\ell_2$-sum
of $\ell_{p_n}^{\sqrt {k_n}}$.  Now for $p>2$, $\ell_p^k$ has a
subspace whose GL-constant is of the same order as the distance
of $\ell_p^k$ to $\ell_2^k$.  Thus if $X$ is not isomorphic to a
Hilbert space, then $W$ has a subspace $V$ which does not have
the GL-property.  But $X\oplus W$ is finitely crudely
representable in $X$.  Thus $X$ is isomorphic to a Hilbert space
if $X$ is determined by its finite dimensional subspaces.
\hfill\qed

\bs

\noindent{\bf 8. Problems}

\ms

\noindent (1) \ Say that a Banach space  {\sl is
determined by its  uniform structure\/} if it is
isomorphic to every Banach space to which it is
uniformly equivalent. It seems that most ``natural" or
``classical" separable spaces are determined by their
uniform structure but at present the tools to establish
this fact are  limited.  Here are some specific
problems, the first three being the most interesting.
\sm
\item{(a)} Is $L_p(0,1)$, $1<p<\infty$, $p\not =
2$,  determined by its uniform structure?
\sm
\item{(b)} Are $c_0$ and $\ell_1$  determined by
their uniform structure?  For $c_0$ partial results are
given in section 3.  For $\ell_1$ we have practically no
information beyond Ribe's theorem [Rib1].
\sm
\item{(c)} Are separable $C(K)$ spaces and $L_1(0,1)$
 determined by their uniform structure?  In
particular, if
$K_1$ and $K_2$ are compact metric spaces for which
$C(K_1)$ is uniformly homeomorphic to $C(K_2)$, must
$C(K_1)$ be isomorphic to $C(K_2)$?
\sm
More generally, one can ask

\item{(d)} Is every separable $\sLp$ space, $1\le p
\le \infty$, determined by its uniform structure?
For the case $p=1$ it is even open whether a uniform
homeomorph of a $\sLo$ space must itself be a $\sLo$
space.  For $1<p<\infty$, $p\not=2$, the simplest
unknown case is the space $\ell_p\oplus \ell_2$.
\sm
In connection with the results of section 2 it is
natural to ask

\item{(e)} Is $\ell_p\oplus\ell_r$ for
$1<p<2<r<\infty$ determined by its uniform structure?

\ms

\noindent  (2) \ A general question related to
uniform homeomorphisms is whether the uniform
structure of a Banach space is determined by the
structure of its discrete subsets. Suppose that $A$
is an $a$-separated $b$-net for a Banach space $X$
with $0<a<b<\infty$; that is, every two points in $A$
are of distance apart larger than $a$ and every point
in
$X$ is of distance less than $b$ from some point of
$A$. If
$T$ is a uniform homeomorphism from $X$ onto $Y$, then
$T[A]$ is a $c$-separated $d$-net for $Y$ for some
$0<c<d<\infty$,  and, since $T$ and $T^{-1}$ are
``Lipschitz for large distances",  the restriction
of $T$ to $A$ is bi-Lipschitz.  We ask whether the
converse is true:

\item{} If there exist $b$-nets $A$ and $B$ in $X$ and
$Y$, respectively, for some $b<\infty$ so that
there is a bi-Lipschitz mapping from $A$ onto $B$,
then must $X$ and $Y$ be uniformly homeomorphic?

An affirmative answer to this problem would of course
imply that ultrapowers of $X$ and $Y$ are
Lipschitz equivalent; this much at least is true:

\proclaim Proposition 8.1. Suppose that $A$ and $B$
are $b$-nets in $X$ and $Y$, respectively, and there
is a bi-Lipschitz mapping $T$ from $X$ onto $Y$.
Let $\U$ be a free ultrafilter on the natural
numbers.  Then $\UX$ is Lipschitz equivalent to
$\UY$.

\Proof Define a mapping $T_n : {1\over n}A\to
{1\over n}B$ by $T_n(x)={1\over n}T(nx)$.  Then for
each $n$, $T_n$ and $T_n^{-1}$ have the same
Lipschitz constants as $T$ and $T^{-1}$,
respectively.  Given a bounded sequence
$\tilde{x}=\xn$ in $X$, we can select a sequence
$\tilde{y}=\yn$ with $y_n$ in  ${1\over n}A$ so that
$\nm{x_n-y_n}<{b\over n}$.  So in the ultrapower
$\UX$, $\tilde{x}=\tilde{y}$.  Since the $T_n$'s are
uniformly Lipschitz, if for $\tilde{x}$ and
$\tilde{y}$ as above we set
$\tilde{T}(\tilde{x})= \{T_n(y_n)\}_{n=1}^\infty$,
then $\tilde{T}$ is a well defined Lipschitz mapping
from $\UX$ to $\UY$.  Doing the the analogous thing
with the $T_n^{-1}$'s, we see that  $\tilde{T}$ is
in fact a bi-Lipschitz mapping from $\UX$ onto $\UY$.
\hfill\qed

\ms

\noindent (3) \ We recall here a well known problem.
Are any two  separable Lipschitz equivalent
Banach spaces isomorphic?

\ms

\noindent (4) \ Besides uniformly continuous and
Lipschitz maps, there are several other interesting
kinds of nonlinear maps between Banach spaces which
appear in the literature.  A notion which arises
naturally from the theory of quasiconformal maps and
which was introduced and studied in the context of
general metric spaces by Tukia and V\"ais\"al\"a
[TV] is that of quasisymmetric maps.  There are
several variants of this notion, but for our
purpose, in the context of Banach spaces, we can use
the following definition.  A map $f$ from a subset
$G$ of a Banach space $X$ into a Banach space $Y$ is
{\sl quasisymmetric} if there is a constant $M$ so
that whenever $u$, $v$, $w$ are in $G$ with
$\nm{u-v}\le\nm{u-w}$, then
$\nm{f(u)-f(v)}\le M\nm{f(u)-f(w)}$.  V\"ais\"al\"a
asked the following:  For $1\le p < q<\infty$, is
there a homeomorphism between $\ell_p$ and $\ell_q$
which is bi-quasisymmetric?  (More generally one can
ask if $\ell_p$ is determined by its
``quasisymmetric structure".)  In this connection we
mention an unpublished observation of V\"ais\"al\"a
and the second author:  Suppose that $G$ is an open
subset of a separable Banach space $X$ and $f$ is a
quasisymmetric Lipschitz map from $G$ into a Banach
space $Y$ which has the Radon-Nikodym property (in
particular, $Y$ can be any separable conjugate
space).  Then, unless $X$ is isomorphic to a
subspace of $Y$, $f$ must be a constant.  Indeed,
if $T$ is the Gateaux derivative of $f$ at some
point, then $T$ is either $0$ or an into
isomorphism. Hence if $X$ is not isomorphic to a
subspace of $Y$, then the Gateaux derivative of $f$
must be $0$ whenever it exists.  The assumptions on
$X$, $Y$, and $f$ guarantee that the Gateaux
derivative exists almost everywhere (e.g. in the
sense that the complement has measure $0$ with
respect to every nondegenerate Gaussian measure on
$X$; see [Ben]).  By Fubini's theorem it follows
that the restriction of $f$ to a line $x+L$ has
derivative $0$ almost everywhere for almost all
translates $x+L$ of $L$.  This implies that $f$ is
constant.

The assumption that $Y$ has the Radon-Nikodym
property is essential, since Aharoni [Aha] proved
that there is a Lipschitz embedding of $C(0,1)$ into
$c_0$.  Of course, every Lipschitz embedding  is
quasisymmetric.

\ms

\noindent (5) \ Theorem 5.8 gives  examples
of spaces whose uniform structure determine exactly
$2^k$ isomorphism classes for $k=0,1,2,\dots$.  From
the construction of [AL2] one easily gets examples
of spaces whose uniform structure determine
$2^{\aleph_0}$ isomorphism classes.  If $\a$ is
any cardinal less than the continuum   which is not
a power of two, in
particular if $\a=3$ or $\a=\aleph_0$, we do not
know how to construct a space which determines
exactly $\a$ isomorphism classes.

\ms

\noindent (6) \ If $X$ is
determined by its pavings (as defined in section 7), must
$X$ be determined by its uniform structure?  More generally,
if two spaces are uniformly homeomorphic, must they have a
common paving?  (In the terminology of section 7,  they do
have the same finite dimensional subspaces by Ribe's
theorem  [Rib1].)  Notice that Corollary 7.7 and the second
remark after Proposition 2.9 give examples of spaces which are
determined both by their pavings and by their uniform
structure.

\bs
\centerline{\bf References}

\medskip

\item{[Aha]} I. Aharoni, {\sl Every separable metric space is Lipschitz
equivalent to a subset of $c_0$,\/} {\bf Israel J. Math. 19} (1974),
284--291.

\item{[AL1]} I. Aharoni and J. Lindenstrauss, {\sl Uniform
equivalence between Banach spaces,} {\bf Bull. Amer. Math. Soc.
84} (1978), 281--283.

\item{[AL2]} I. Aharoni and J. Lindenstrauss, {\sl An extension of
a result of Ribe,\/} {\bf Israel J. Math. 52} nos. 1-2 (1985),
59--64.

\item{[Ald]} D. Aldous, {\sl Subspaces of $L^1$ via random
measure,\/} {\bf Trans. Amer. Math. Soc. 258} (1981), 445--463.

\item{[Als]} D. E. Alspach, {\sl Quotients of $C[0,1]$ with separable
dual,\/} {\bf Israel J. Math. 29} no. 4 (1978), 361--384.

\item{[AB]} D. E. Alspach and Y. Benyamini, {\sl $C(K)$ quotients of
separable ${\cal L}_\infty$ spaces,\/} {\bf Israel J. Math. 32} nos. 2-3
(1979), 145--160.

\item{[Bak]} J. W. Baker, {\sl Dispersed
images of topological spaces and uncomplemented subspaces of $C(X)$,\/}
{\bf Proc. Amer. Math. Soc. 41} (1973), 309--314.

\item{[BDGJN]}  G. Bennett, L.~E. Dor, V. Goodman, W.~B. Johnson, and C.~M.
Newman, {\sl  On uncomplemented subspaces of  $L_p$, $1 < p < 2$,} {\bf
Israel J. Math. 26} (1977), 178--187.

\item{[Ben1]} Y. Benyamini, {\sl An extension theorem for separable Banach
spaces,\/} {\bf Israel J. Math. 29} (1978), 24--30.

\item{[Ben2]} Y. Benyamini, {\sl The uniform classification of
Banach spaces,} {\bf Longhorn Notes, Texas Functional Analysis Seminar,
1984-85} Univ. of Texas,  15--39 (An electronic version can be
obtained via FTP from ftp.math.okstate.edu /pub/banach.)

\item{[BL]} J. Bergh and J. L\"ofstr\"om, {\sl Interpolation spaces. An
introduction,\/} Springer-Verlag {\bf Grundlehren Math. Wissenschaften 223}
(1976).

\item{[Bou]} J. Bourgain, {\sl Remarks on the extensions of Lipschitz maps
defined on discrete sets and  uniform homeomorphisms,} {\bf Lecture Notes
in Math. 1267} (1987), 157--167.

\item{[BP]} C. Bessaga and A. Pe\l czy\'nski, {\sl Spaces of continuous
functions (IV),\/} {\bf Studia Math. 19} (1960), 53--62.

\item{[Cal]} A. Calderon, {\sl Intermediate spaces and interpolation, the
complex method,\/} {\bf Studia Math. 24} (1964), 113--190.

\item{[CJT]} P.~G. Casazza, W.~B. Johnson, and L. Tzafriri, {\sl On
Tsirelson's space,} {\bf Israel J. Math. 47} (1984), 81--98.

\item{[CO]} P.~G. Casazza and E. Odell, {\sl Tsirelson's space and minimal
subspaces,\/} {\bf Longhorn Notes, Texas Functional Analysis Seminar,
1982-83} Univ. of Texas, 61--72.

\item{[CS]} P.~G. Casazza and T. J. Shura, {\sl Tsirelson's space,}
{\bf Lecture Notes in Math. 1363} Springer-Verlag  (1989)

\item{[DGZ]} R. Deville, G. Godefroy, and V.~E. Zizler, {\sl The three
space problem for smooth partitions of unity and $C(K)$ spaces,\/} {\bf
Math. Ann. 288} (1990), 613--625.

\item{[Ede]} I. S. Edelstein, {\sl On complemented subspaces and
unconditional bases in $\l_p+\l_2$,\/} {\bf Teor. Funkcii, Functional.
Anal. i Prilozen. 10} (1970), 132--143 (Russian).

\item{[EW]} I. S. Edelstein and P. Wojtaszczyk, {\sl On projections and
unconditional bases in direct sums of Banach spaces,} {\bf Studia Math.
56} (1976), 263--276.

\item{[Enf]} P. Enflo, {\sl On the non-existence of uniform
homeomorphisms   between  $L_p$ spaces,} {\bf  Ark. Mat. 8} (1969),
195--197.

\item{[FJ]} T. Figiel and W. B. Johnson, {\sl A uniformly convex Banach
space which contains no $\ell_p$,}\/  {\bf Compositio Math. 29} (1974),
179--190.

\item{[Gor]} E. Gorelik, {\sl The uniform nonequivalence of
$L_p$ and $\ell_p$,} {\bf Israel J. Math. 87} (1994), 1--8.

\item{[HM]} S. Heinrich and P. Mankiewicz, {\sl Applications
of ultrapowers to the uniform and Lipschitz classification of
Banach spaces,} {\bf Studia Math. 73} (1982), 225--251.

\item{[Hei]}  S. Heinrich, {\sl Ultraproducts in Banach space theory,}
{\bf J. f\"ur Die Reine und Angewandte Math. 313} (1980), 72--104.

\item{[Jam]} R.~C.~James, {\sl Some self-dual properties of normed linear
spaces,\/} {\bf Annals  Math. Studies 69} (1972).

\item{[Joh1]} W. B. Johnson, {\sl Operators into  $L_p$ which factor
through    $\ell_p$,} {\bf J. London Math. Soc. 14} (1976), 333--339.

\item{[Joh2]}  W. B. Johnson, {\sl A reflexive Banach space which is not
sufficiently Euclidean,}\/  {\bf Studia Math. 55} (1976), 201--205.

\item{[Joh3]}  W. B. Johnson,{\sl Banach spaces all of whose subspaces have
the approximation property,} {\bf  Special Topics of Applied  Mathematics}
North-Holland (1980), 15--26.

\item{[JO]} W. B. Johnson and E. Odell, {\sl Subspaces of
$L_p$ which embed into $\ell_p$,} {\bf Compositio Math. 28} (1974),
37--49.

\item{[JSc]}   W. B. Johnson and G. Schechtman, {\sl On the distance of
subspaces of
$l^n_p$ to $l^k_p$,}  {\bf Trans.  Amer. Math. Soc. 324} (1991), 319--329.

\item{[JS]} W. B. Johnson and A. Szankowski, {\sl Complementably universal
Banach spaces,\/} {\bf Studia Math. 58} (1976), 91--97.

\item{[JZ]} W. B. Johnson and M. Zippin, {\sl On subspaces of quotients of
$(\Sigma G)_{\ell _p}$ and $(\Sigma G) _{c_0}$,}\/  {\bf Israel J. Math. 13
nos. 3 and 4} (1972), 311--316.

\item{[Kad]} M. I. Kadec, {\sl A proof of the topological
equivalence of separable Banach spaces,\/} {\bf Funkcional Anal.
i Prilo\v zen 1} (1967), 53--62 (Russian).

\item{[KP]} M. I. Kadec and A. Pe\l czy\'nski, {\sl Bases,
lacunary   sequences and complemented subspaces in the spaces  $L_p$,} {\bf
Studia   Math. 21} (1962), 161--176.

\item{[KT]} R. Komorowski and N. Tomczak-Jaegermann, {\sl   Banach spaces
without local unconditional structure,\/}  {\bf  Israel J. Math. 89}
(1995), 205--226.

\item{[Kri]} J. L. Krivine, {\sl Sous espaces de dimension finite des
espaces de Banach reticul\'es,\/} {\bf Annals Math. 104} (1976), 1--29.

\item{[KM]}  J. L. Krivine and B. Maurey, {\sl Espaces de Banach
stables,\/} {\bf Israel J. Math. 39} (1981), 273--295.

\item{[Lin1]} J. Lindenstrauss, {\sl On non-linear projections in Banach
spaces,} {\bf Mich. J. Math. 11} (1964), 268--287.

\item{[Lin2]} J. Lindenstrauss, {\sl
On non-separable reflexive Banach spaces,} {\bf Bull. Amer. Math. Soc. 72}
(1966), 967--970.

\item{[LPe]} J. Lindenstrauss and A. Pe\l czy\'nski, {\sl Absolutely summing
operators in ${\cal L}_p$ spaces and their applications,} {\bf Studia Math.
29} (1968), 275--326.

\item{[LP]} J. Lindenstrauss and D. Preiss, {\sl On almost Frech\'et
differentiability of Lipschitz functions,}

\item{[LR]} J. Lindenstrauss and H.~P. Rosenthal, {\sl The  ${\cal L}_p$
spaces,} {\bf Israel J. Math. 7} (1969), 325--349.

\item{[LT1]} J. Lindenstrauss and L. Tzafriri, {\sl Classical Banach spaces
I,   Sequence spaces,}  Springer-Verlag, (1977).

\item{[LT2]} J. Lindenstrauss and L. Tzafriri, {\sl Classical Banach spaces
II,   Function spaces,}  Spring\-er-Verlag, (1979).

\item{[Loz]} G. Ya. Lozanovskii, {\sl On some Banach lattices,\/} {\bf
Siberian Math. J. 10} (1969), 419--431 (English translation).

\item{[Mau1]} B. Maurey, {\sl Th\'eor\`emes de
factorisation pour les op\'erateurs lin\'eaires \`a valeurs dans les
espaces $L^p$} {\bf Ast\'erisque 11} (1974).

\item{[Mau2]} B. Maurey, {\sl  Type et cotype
dans les espaces munis de structures locales inconditionnelles,\/} {\bf
Seminaire   Maurey-Schwartz 1973-74}  Expose 24-25, Ecole Polytechnique,
Paris.

\item{[Mau3]} B. Maurey, {\sl  Un theor\`eme de prolongement,}
{\bf C.~R.~Acad., Paris 279} (1974), 329--332.

\item{[MP]} B. Maurey and G. Pisier, {\sl  S\'eries de variables
aleatoires vectorielles independantes et propri\'et\'es geometriques des
espaces de Banach,\/} {\bf Studia Math. 58} (1976), 45--90.

\item{[Maz]} S. Mazur, {\sl Une remarque sur l'hom\'eomorphismie des champs
fonctionnels,\/} {\bf Studia Math. 1} (1930), 83--85.

\item{[Mic]} E. Michael, {\sl Continuous selections, I,\/} {\bf Annals Math.
63} no. 2 (1956), 361--382.

\item{[OS]} E.~W. Odell and T. Schlumprecht, {\sl The distortion
problem,\/} {\bf Acta Math.}  (to appear).

\item{[PS]} A. Pe\l czy\'nski and Z. Semadeni, {\sl Spaces of continuous
functions, III,\/} {\bf Studia Math. 25} (1959), 211--222.

\item{[Pis]}  G. Pisier, {\sl The volume of convex bodies and Banach space
geometry,}  {  Cambridge University Press} (1989).

\item{[Rib1]} M. Ribe, {\sl On uniformly homeomorphic normed
spaces,} {\bf Ark. Math. 14} (1976), 237--244.

\item{[Rib2]} M. Ribe, {\sl On uniformly homeomorphic normed
spaces, II,} {\bf Ark. Math. 1} (1978), 1--9.

\item{[Rib3]} M. Ribe, {\sl Existence of separable uniformly homeomorphic
non isomorphic Banach spaces,} {\bf Israel J. Math. 4} (1984), 139--147.

\item{[Sem]} Z. Semadeni, {\sl Banach spaces of continuous functions,\/}
Polish Scientific Publishers, Warsaw (1971).

\item{[T-J]} N. Tomczak-Jaegermann, {\sl Banach-Mazur distances and
finite-dimensional operator ideals,} {\bf Pitman Monographs and Surveys in Pure
and Applied Mathematics  38}, Longman (1989).

\item{[Tor]} H. Torunczyk, {\sl Characterizing Hilbert space topology,\/}
{\bf Fund. Math. 111} (1981), 247--262.

\item{[Tsi]} B. S. Tsirelson, {\sl Not every Banach space contains an
embedding of $\ell_p$ or $c_0$,\/} {\bf Functional Anal. Appl. 8} (1974),
138--141 (translated from the Russian).

\item{[TV]} P. Tukia and J. V\"ais\"al\"a, {\sl Quasisymmetric embeddings
of metric spaces,\/} {\bf Ann. Acad. Sci. Fenn. Ser. A I Math. 5} (1980),
97--114.

\item{[Woj]} P. Wojtaszczyk, {\sl On complemented subspaces and
unconditional bases in $\l_p+\l_q$,\/} {\bf Studia Math. 48} (1973),
197--206.

\item{[Zip]} M. Zippin, {\sl  The separable extension problem,\/} {\bf
Israel J. Math. 26} nos. 3-4 (1977), 372--387.

\ms

\settabs\+College Station, TX  77843--3368 U.S.A.aaaa&The
Hebrew University of Jerusalem\cr

\+W. B. Johnson&J. Lindenstrauss\cr
\+Department of Mathematics&Institute of Mathematics\cr
\+Texas A\&M University&The
Hebrew University of Jerusalem\cr
\+College Station, TX  77843--3368 U.S.A.&Jerusalem, Israel\cr
\+{\tt
email: johnson@math.tamu.edu}&{\tt email: joram@math.huji.ac.il}\cr
\ms
\line{G. Schechtman\hfill}
\line{Department of Theoretical Mathematics\hfill}
\line{The Weizmann Institute of Science \hfill}
\line{Rehovot, Israel\hfill}
\line{{\tt email: mtschech@weizmann.weizmann.ac.il}\hfill}

\bye